\documentclass[12pt]{article}
\usepackage{latexsym,amssymb,amsmath,amsfonts,amsthm}
\usepackage{color}

\newtheorem{theorem}{Theorem}

\newtheorem{lemma}{Lemma}

\newtheorem{corollary}{Corollary}

\def\beq{ \begin{equation} }
\def\eeq{ \end{equation} }
\def\mn{\medskip\noindent}
\def\ms{\medskip}

\def\ep{\epsilon}

\def\square{\vcenter{\vbox{\hrule height .4pt
  \hbox{\vrule width .4pt height 5pt \kern 5pt
        \vrule width .4pt} \hrule height .4pt}}}

\def\sqz{\kern -0.2em}

\def\E{\mathbb{E}}
\def\P{\mathbb{P}}
\def\RR{\mathbb{R}}
\def\TT{\mathbb{T}}
\def\ZZ{\mathbb{Z}}

\def\clearp{}

\def\iint{\int\kern -0.5em \int}

\definecolor{darkgreen}{rgb}{0,0.7,0}

\usepackage{graphicx}
\graphicspath{%
    {converted_graphics/}
    {/}
}
\begin{document}

\title{The q-voter model on the torus}
\author{Pooja Agarwal, Mackenzie Simper, and Rick Durrett 
\\
Brown, Stanford, and Duke}
\date{\today}						

\maketitle

\begin{abstract}
In the $q$-voter model, the voter at $x$ changes its opinion at rate $f_x^q$, where $f_x$ is the fraction of neighbors with the opposite opinion. Mean-field calculations suggest that there should be coexistence between opinions if $q<1$ and clustering if $q>1$. This model has been extensively studied by physicists, but we do not know of any rigorous results. In this paper, we use the machinery of voter model perturbations to show that the conjectured behavior holds for $q$ close to 1. More precisely, we show that if $q<1$, then for any $m<\infty$ the process on the three-dimensional torus with $n$ points survives for time $n^m$, and after an initial transient phase has a density that it is always close to 1/2. If $q>1$, then the process rapidly reaches fixation on one opinion. It is interesting to note that in the second case the limiting ODE (on its sped up time scale) reaches 0 at time $\log n$ but the stochastic process on the same time scale dies out at time $(1/3)\log n$. 
\end{abstract}

\section{Introduction}

In the linear voter model, the state at time $t$ is $\xi_t :\ZZ^d \to \{0,1\}$, where 0 and 1 are two opinions. The individual at $x$ changes opinion at a rate equal to the fraction $f_x$ of its neighbors with the opposite opinion. For the last decade physicists have studied the $q$-voter model, in which the flip rate at $x$ is $f_x^q$. When $q$ is an integer, the dynamics may be thought of as: select $q$ neighbors of $x$ uniformly, and change the opinion of $x$ if all $q$ neighbors disagree with $x$. However, there is no reason to restrict $q$ to be an integer. Abrams and Strogatz \cite{AbrStr} introduced this system in 2003 as a model of language death, and argued based on data on languages in 42 regions that $q = 1.31 \pm 0.25$. In the physics literature there have been many studies of the system on lattices, complex networks, and even on graphs that co-evolve with the state of individuals. See \cite{CMPS, HCDM, AJ, MSM, MLCPS,  VCSM, VazLop} and references therein. According to \cite{MSM}, for finite but large systems, the process with $q<1$ can remain in a dynamically active phase for observation times that grow exponentially with $n$, while for $q>1$ the transition into an absorbing state is `abrupt'. 

The difference between $q<1$ and $q>1$ is due to the different types of frequency dependence in the two models. When $q<1$, rare opinions spread more rapidly compared to the voter model, while for $q>1$, they spread more slowly. A more quantitative viewpoint is provided by mean field theory. This analysis is often done by writing an equation by pretending sites are always independent of each other. Here, we will instead consider the system on the complete graph in which each site interacts equally with all the others. In this case, the frequency of 1's, $u$, satisfies
$$ 
du/dt = -u (1-u)^q + (1-u) u^q = u(1-u) g(u)
$$ 
where $g(u) = u^{q-1} - (1-u)^{q-1}$. This system has three fixed points: $0$, $1/2$ and $1$. 

\begin{itemize}
  \item 
If $q < 1$, $g(u)$ decreases from $\infty$ to $-\infty$ as $u$ increases from 0 to 1. So the fixed points 0 and 1 are unstable and the interior one is attracting. In this case it is expected that coexistence occurs. 

\item 
If $q > 1$, $g(u)$ increases from $-1$ to $1$ as $u$ increases from 0 to 1. So the fixed points 0 and 1 are stable and the interior one is unstable. In this case it is expected that clustering occurs. That is, we will see larger and large regions occupied by one type.
\end{itemize}

\noindent
For more on the heuristics that lead to these conclusions, see the 1994 paper by Durrett and Levin \cite{DL}. In most of the papers in the physics literature, the analysis is done by using the pair-approximation, which is equivalent to supposing that the state of the system is always a Markov chain. 

Recently, Vasconclos, Levin, and Pinheiro \cite{VLP} have considered a version of the $q$-voter in which the powers $q_1$ and $q_0$ for flipping to 1 and 0 can be different. They did this to study complex contagions which have been used to model the spread of idioms and hashtags on Twitter \cite{RMK} and in many other situations, see the book by Centola \cite{HBS}. When $q_1 \neq q_0$, there arises situations when one opinion dominates the other, see Figure 2a in \cite{VLP}, but the situation with $q_1=q_0$ seems to capture of all of the interesting behavior.

\subsection{Voter model perturbations}

The linear voter model has a rich theory due to its duality with coalescing random walk. This duality exists because the process can be constructed from a graphical representation. See Section \ref{ss:vmd} for details.  However, the inherent asymmetry between 1's and 0's in the graphical representation makes it impossible to construct nonlinear voter models where the flip rates depend only on $f_x$. See Section \ref{ss:nlv} for a proof. 

To get around this difficulty, we will suppose $q$ is close to 1 and view the system as a voter model perturbation in the sense of Cox, Durrett, and Perkins \cite{CDP}. On $\ZZ^d$, this theory requires $d\ge 3$ so that the voter model has a one parameter family of stationary distributions $\nu_u$, $0 \le u \le 1$. For this and other elementary facts about the voter model that we use, see Liggett's 1999 book \cite{L99}. 

In general, the rate of flipping from $i$ to $j\neq i$ in a voter perturbation has the form
$$
c^\delta_{i,j}(x,\xi) = f_j + \delta^2 h_{i,j}(x,\xi)
$$
where $f_j$ is the fraction of neighbors in state $j$, and $h_{i,j}(x,\xi)$ is the perturbation to the rate of flipping from $i$ to $j$. Usually the perturbation variable is $\ep$, but here it will be convenient to let $\ep=\delta^2$. To simplify formulas we will assume $h_{i,j}(x,\xi)=0$ when $\xi(x)\neq i$. Here we will consider the special case in which the neighborhood has size $k$ and the flip rate only depends on the number of neighbors $n(x)$ in state $j$:
$$
c^\delta_{i,j}(x,\xi) = f_j + \delta^2 r^k_{n(x)} \qquad\hbox{for $1\le n(x) \le k$}.
$$
The $r^k_i$ do not have to be nonnegative, see (1.7) in \cite{CDP}, but we will suppose $r^k_0=0$ so that $\equiv 0$ and $\equiv 1$ are absorbing states. For simplicity, we will restrict our attention to three dimensions. In that context, we will consider neighborhoods $x +{\cal N}$ with $0 \notin {\cal N}$ and $|{\cal N}| \ge 3$ chosen so that the group generated by ${\cal N}$ is $\ZZ^3$

\mn
{\bf q-voter model.} The rate at which a site $x$ flips to 0 in the $q$-voter model is $f_x^q$, where $f_x$ is the fraction of neighbors with the opposite opinion. Suppose for the moment that $q<1$. In this case, if we write
$$
f_x^q = f_x + (f_x^q - f_x),
$$
then the term in parentheses is $\ge 0$. Let $q=1 -\delta^2$ and write $u$ instead of $f_x$ Then, 
$$
 u^q - u =  u \left( u^{-\delta^2} - 1 \right) = u  \left( \exp(\delta^2 \log(1/u)) - 1 \right) \approx \delta^2 u \log(1/u).
$$
From this we see that if $q<1$, then the perturbation is 
\beq
r^k_i = (i/k) \log(k/i).
\label{rforqv}
\eeq
which vanishes when $i=0$ or $k$. 

If we let $q=1 + \delta^2$ and again write $u$ instead of $f_x$, then
$$
 u^q - u =  u \left( u^{\delta^2} - 1 \right) = u  \left( \exp(\delta^2 \log(u)) - 1 \right) \approx - \delta^2 u \log(1/u).
$$
Hence when $q>1$, the perturbation is
\beq
r^k_i = -(i/k) \log(k/i).
\label{rforq>1}
\eeq

\subsection{ODE limit} \label{ss:ODElim}

Following the approach of Cox and Durrett \cite{CD}, who used the voter perturbation machinery to study evolutionary games on the torus in dimension $d \ge 3$, we  will consider the $q$-voter model in what they called the weak-selection regime. (For results in the strong selection regime see Section \ref{ss:strsel}.) Let $\TT_n$ be the three dimensional torus with $n$ points and hence side length $L=n^{1/3}$. Let $\ep_n = \delta^2_n$. The first thing to do is to prove convergence of the density of 1's,  
$$
U_n(t) = \frac{1}{n} \sum_{x \in \TT_n} \xi_{t/\ep_n}(x),
$$
to the solution of an ODE. Let $\rho_m^i$ denote the probability that in $\nu_u$ the origin is in state $i$ while exactly $m$ of the neighbors are in state $1-i$.
We write $a_n \ll b_n$ for positive quantities $a_n$ and $b_n$ to indicate $a_n/b_n \to 0$ as $n\to\infty$.

\begin{theorem}\label{detlim}
Suppose $q=1-\ep_n$ with $n^{-1} \ll \ep_n \ll n^{-2/3}$. If $U_n(0) \to u_0$ then $U_n(t)$ converges uniformly on compact sets to the solution of the ODE
\beq
\frac{du}{dt} = \sum_{m=1}^{k-1} r^k_i (\rho^0_m(u) - \rho^1_m(u)) \qquad u(0)=u_0
\label{ODElim}
\eeq
\end{theorem} 

Intuitively, Theorem \ref{detlim} holds due to a separation of time scales. The voter model runs at a fast rate, so when the density is $u$ on the torus, the system has distribution $\approx\nu_u$. The rate of change of the density can then be computed by looking at the expected rate of change when the state is $\nu_u$. Writing $\langle \ \rangle_u$ for expected value with respect to $\nu_u$, the right hand side of the ODE is
\beq
\phi(u) = \langle h_{0,1}- h_{1,0} \rangle _u =  \sum_{m=1}^{k-1} r^k_i (\rho^0_m(u) - \rho^1_m(u)).
\label{ODErhs}
\eeq
This result will be proved by constructing the process on a graphical representation and then defining a dual that is a coalescing branching random walk. The voter part of the process leads to a coalescing random walk. When a perturbation event occurs at a point $x$, the dual branches to include all of the points in $x+{\cal N}$.  This will be described in detail in Section \ref{ss:vpdual}. The proof of Theorem \ref{detlim} is almost identical to the proof of Theorem 6 in Cox and Durrett \cite{CD} so we will only outline the proof, referring to \cite{CD} for details. 
When $\ep_n \gg n^{-2/3}$ the particles in the dual have time to wrap around the torus and come to equilibrium in between branching events. It is known that on the torus if we start two random walks from independent randomly chosen locations, then the time to coalesce is of order $n$. Thus the assumption $\ep_n \ll n^{-1}$ is needed for the perturbation to have an effect. 

Computing the $r^i_m(u)$, see Section \ref{sec:cpert}, leads to the following ODE

\begin{theorem} \label{limitODE}
In the three dimensions when the neighborhood has size $k$, the limiting ODE is 
$$
\frac{du}{dt} = \pm c_k u(1-u)(1-2u)f_k(u)
$$
where $f_k(u)$ is a polynomial that is positive on $[0,1]$ and $f(0)=f(1)=1$.
We have $+$ for $q<1$ and $-$ for $q<1$.
\end{theorem}

When $q<1$, the fixed point at 1/2 is attracting and we have

\begin{theorem}\label{persist}
Suppose $q=1-\ep_n$ and $\ep_n \sim Cn^{-a}$ for some $a\in(2/3,1)$. There is a $T_0$ that only depends on $u_0$, so that for any $\gamma>0$ and $m<\infty$, if $n$ is large then with high probability
$$
| U_n(t) - 1/2 | \le \gamma  \hbox{ for all $t\in[T_0,n^m]$}.
$$
\end{theorem}

\noindent
Here and in what follows ``with high probability'' means with probability $\to 1$ as $n\to\infty$. 

To prove Theorem \ref{persist}, we will follow the approach of Huo and Durrett \cite{latent} who proved a similar result for the latent voter model on a random graph generated by the configuration model. Although the random graph has a more complicated geometry than the torus, the proof in that setting is simpler than the one given here, since on the graph random walks mix in time $O(\log n)$ rather that in time $O(n^{2/3})$.

\clearpage

\begin{figure}[tbp] 
  \centering
  \includegraphics[width=3.92in,height=3.92in,keepaspectratio]{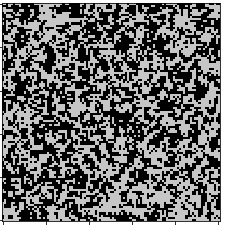}
  \caption{Cross-section from a simulation of $q=0.9$ on a $100 \times 100 \times 100$ grid with periodic boundary conditions.}
  \label{fig:q09}
\end{figure}

\ms
{\bf Outline of the proof of Theorem \ref{persist}.} 

\begin{itemize} 
  \item 
Section \ref{ss:DarNor} introduces a general result for proving convergence of stochastic processes to limiting ODEs, due to Darling and Norris \cite{DN}, which is the key to the proofs of the persistence results for our model (and for the latent voter model). The main difficulty is to bound the difference between the drift in the density $U_n$ of the particle system and the drift in the ODE. In particular, one must prove that the drift in the density of $U_n$, which is a function of the configuration, is almost a function of the overall density.

\item
In Section \ref{ss:igbr} we take the first step in the proof, which is to show that if $2/3 < b < a$ then we can ignore the perturbation on  $[t/\ep_n-n^b,t/\ep_n]$, i.e., the process will evolve like the voter model. This has the consequence that if there are $n \cdot u$ 1's at time $t/\ep_n - n^b$, then at time $t/\ep_n$ the process is close to the voter equilibrium $\nu_u$. The argument here is an improvement over the one in Section 3.1 of \cite{latent}. We use Azuma's inequality to get error estimates that are stretched exponentially small, i.e., $\le C \exp(cn^{-\alpha})$ with $\alpha>0$ rather than polynomial, i.e., $\le Ct^{-p}$. 

\item
In Section \ref{ss:dens} we introduce a result about ``renormalizing'' the voter model, that comes from work of Bramson and Griffeath \cite{BGrenorm} in $d=3$ and Z\"ahle \cite{Zrenorm} in $d \ge 3$. They show that if we consider the number of 1's in the voter model equilibrium with density $\lambda$, $\xi^\lambda$, in a cube $Q(r)$ of side $r$, then 
\beq
\widehat S_r = (\lambda(1-\lambda))^{-1/2} r^{-5/2} \left(\sum_{x \in Q(r)} \xi^\lambda(x) - \lambda\right) \Rightarrow \hbox{Normal}(0,C)
\label{vlamclt}
\eeq
We use this to obtain information about a similar normalized sum $T_r$ of the number of ones in a cube of side $r$ on the torus at time $t/\ep_n$ when the number of 1's at time $t/\ep_n -n^b$ is $\lambda n$. To be specific, we let $\bar S_n$ be the normalized sum of $\xi^\lambda_{\sigma(n)}(x)$ in the process that starts at time 0 from product measure with density $\lambda$ and is run for time $\sigma(n) = n^{0.6}$. We show that $\bar S_r \le T'_r \le \widehat S_r$, where $T'_r$ is a small modification of $T_r$. 

\item
In Section \ref{ss:TvsT} we bound the difference between $T'_r$ and $T_r$. This in turn gives us a bound on the largest coalescing random walk cluster in $T_r$ in $Q(r)$, see \eqref{maxZbd}, and a bound on the fluctuations of the density in the cubes, which is important for completing the next step.

\item
In Section \ref{ss:diffm} we bound the difference between the drifts in the particle system and the ODE. To do this, we have to show that the empirical finite distributions on the torus $\TT_n$ are close to the values that come from $\nu_u$. In doing this we rely on the result about the density in cubes proved in Section \ref{ss:dens} to divide space at time $t/\ep_n + s_n$ into cubes with $n^{b(3)}$ sites, where $b(3) > b(2)$.  Here $s_n = n^{(2+\alpha)b(2)/3}$ with $\alpha$ small, so that the empirical f.d.d.'s in cubes of volume $n^{b(3)}$ that do not touch are almost independent. This leads to errors of size $C\exp( - n^{1-b(3)-2\alpha})$. 

\item
In Section \ref{ss:fd} we put the pieces together to prove the result. As in Section 3.5 of \cite{latent} we do this by showing that if the density $U_t$ reaches $|U_t -1/2| = 4\ep$ then with very high probability (i.e., for any $k$ with probability $\ge 1-n^{-k}$ for large $n$) it will return to $|U_t -1/2| \le \ep$ before we have $|U_t -1/2|> 5\ep$. Taking $\delta=5\ep$ gives the desired result

\end{itemize}

In all of our estimates except those in Sections \ref{ss:dens} and \ref{ss:TvsT}, the errors are bounded stretched exponentially small, so we 

\mn
{\bf Conjecture.} {\it When $q<1$ the process persists for time $\exp(n^\beta)$ for some $\beta>0$.}

\mn
The could be proved with a rather small value of $\beta$ if the errors in \eqref{mombdT} and  \eqref{maxZbd} could be improved to be stretched exponentially small.
Readers familiar with long time survival results for the contact process, see e.g., Section 3 in part I of Liggett \cite{L99}, might expect the conjecture to say survival occurs for time $\exp(\gamma n)$ with $\gamma>0$. However, the conjecture above cannot hold for $\beta > 1/3$. If we run time backwards from $t/\ep_n$ to $t/\ep_n - n^{2/3}$ then the $n$ initial particles in the CRW will have coalesced to $n^{1/3}$ particles. If all of these happen to land on sites in state 0 at time $t/\ep_n - n^{2/3}$ the process will go extinct at time $t/\ep_n$.

\subsection{Rapid Extinction when $q>1$} \label{ss:rapid}

When $q>1$, the fixed point at 1/2 is unstable while the ones at 0 and 1 are locally attracting
To get rid of the constant $c_k$ in the ODE limit we consider 
$$
U_n(t) = \frac{1}{n} \sum_{x \in \TT_n} \xi_{t/\ep_n c_k}(x)
$$

\begin{theorem}\label{dieout}
Suppose $q=1+\ep_n$ and $\ep_n \sim Cn^{-a}$ for some $a\in(2/3,1)$. If $U_n(0) = u_0 < 1/2$ and $\alpha> 1/3$ then
$$
P \left( U_n(\alpha \log n)=0 \right) \to 1 \qquad\hbox{as $n\to\infty$.}
$$ 
\end{theorem}

\noindent
This is proved in Section 5. Much of the work for the proof of Theorem \ref{dieout} has already been done in the proof of Theorem \ref{persist}. Those results imply that the density in the particle system stays close to the solution of the ODE. To be precise, we can show that with high probability. 
$$
|U_n(t) - u(t)| \le \ep u(t) \quad\hbox{until}\quad \tau = \inf\{ t : x_t \le n^{-(1-b(0))} \}
$$
where $2/3 < b(0) < \min\{b,1-\alpha\}$. Since the ODE is $u'(t) = - f(u)$ with $f(u)/u \to 1$ as $u \to 0$, the limiting ODE has $u(\alpha \log n) \approx n^{-\alpha}$. Our proof shows that when the density gets to $\le n^{-b(0)}$ fluctuations in the voter model make the system go extinct in a time that is $\le Cn^b$. See Section \ref{sec:rapext} for details. The keys to the voter extinction result are (i) the observation that the number of 1's in the voter model is a time change of continuous-time symmetric random walk, and (ii) results on the size of the boundary of the voter model in the low density regime due to Cox, Durrett, and Perkins  \cite{CDP2}.

\begin{figure}[tbp] 
  \centering
 \includegraphics[width=4in,height=4.02in,keepaspectratio]{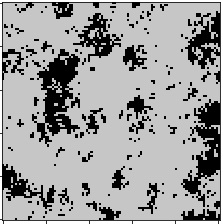}
  \caption{Simulation of $q=1.1$ on a $100 \times 100 \times 100$ grid.}
  \label{fig:q11}
\end{figure}

\subsection{Results for strong selection} \label{ss:strsel}

Let $\xi^\ep_t$ be a voter model perturbation on $\ZZ^d$ with flip rates
$$
c^{\delta_n}_{i,j} = f_j + \delta^2_n h_{i,j}(x,\xi)
$$
where $f_j$ is the fraction of neighbors in state $j$ and the second term is the perturbation. As before we let $\ep_n = \delta_n^2$. In this section we will examine the case $\ep_n \gg n^{-2/3}$, which we call the strong selection regime.

Intuitively, the next result says that if we rescale space to $\delta_n\TT_n$ (recall $\TT_n$ is the three dimensional torus) and speed up time by $\delta_n^{-2}$, then the process converges to the solution of a partial differential equation on $\RR^3$. The torus turns into $\RR^3$ in the limit because $\delta_n \ll n^{-1/3}$ while the torus has side $n^{1/3}$. To make a precise statement, the first thing we have to do is to define the mode of convergence. To simplify the writing we drop the subscript $n$ on $\delta$. Given $r\in(0,1)$, let $a_\delta = \lceil \delta^{r-1} \rceil \delta$, $Q_\delta = [0, a_\delta)^3$, and
$|Q_\delta|$ the number of points in $Q_\delta$. For $x \in a_\delta \ZZ^d$ and $\xi\in \Omega_\delta$, the space of all functions
from $\delta\ZZ^3$ to $S$, let
$$
D_i(x,\xi) = |\{ y \in Q_\delta : \xi(x+y) = i \}|/|Q_\delta|.
$$ 

We endow $\Omega_\delta$ with the $\sigma$-field ${\cal F}_\delta$ generated by the finite-dimensional distributions. Given a sequence of measures
$\lambda_\delta$ on $(\Omega_\delta,{\cal F}_\delta)$ and continuous functions $w_i$, we say that $\lambda_\delta$ has asymptotic densities
$w_i$ if for all $0 < \eta, R < \infty$ and all $i\in S$
$$
\lim_{\delta\to 0} \sup_{x\in a_\delta\ZZ^3, |x| \le R} \lambda_\delta ( |D_i(x,\xi) - w_i(x)| > \eta ) \to 0.
$$

\begin{theorem} \label{hydro}
Suppose $d = 3$.
Let $w_i: \RR^d \to [0,1]$ be continuous with $\sum_{i\in S} w_i = 1$. Suppose the initial conditions $\xi^\delta_0$
have laws $\lambda_\delta$ with asymptotic densities $w_i$ and let  
$$
u^\delta_i(t,x) = P( \xi^\delta_{t\delta^{-2}}(x) = i)
$$ 
If $x_\delta \to x$ then $u^\delta_i(t,x_\delta) \to u_i(t,x)$ the solution of the system of partial differential equations:
\beq
\frac{\partial}{\partial t} u_i(t,x) = \frac{\sigma^2}{2} \Delta u_i(t,x) +  \phi_i(u(t,x))
\label{PDElimit}
\eeq
with initial condition $u_i(0,x) = w_i(x)$. The reaction term 
\beq
\phi_i(u) = \sum_{j \neq i} \langle h_{j,i}(0,\xi) -  h_{i,j}(0,\xi) \rangle_u
\label{phidef}
\eeq
where the brackets are expected value with respect to the voter model stationary distribution $\nu_u$
in which the densities are given by the vector $u$.
\end{theorem}

\noindent
This result is Theorem 2 in \cite{CD}. For more details see that paper.

The intuition is similar to that for the ODE limit in Theorem \ref{detlim}.
On the fast time scale the voter model runs at rate $\delta^{-2}$ versus the perturbation at rate 1, so the states of sites near $x$ at time $t$ 
is always close to the voter equilibrium $\nu_{u(t,x)}$. Thus, we can compute the rate of change of $u_i(t,x)$ by assuming the nearby
sites are distributed according to the voter model equilibrium $\nu_{u(t,x)}$.

Cox and Durrett considered evolutionary games on the torus in $d \ge 3$ with game matrix ${\bf 1} + w G$, where ${\bf 1}$ is a matrix of 1's. Their $w$ corresponds to our $\ep_n$.  When $w=0$ the system reduces to the voter model. They found convergence to an ODE when $n^{-1} \ll w \ll n^{-2/d}$ and convergence to a PDE when $w \gg n^{-2/d}$. Their results can be used  prove a PDE limit for our system when $\ep_n \gg n^{-2/d}$. Since there are only two opinions we only need one variable $u_1$, which corresponds to our $u$. The $\phi$ in \eqref{phidef} is the same as the right hand side of our ODE, which should be clear from \eqref{ODErhs}.

In the case of a $2 \times 2$ game with a stable mixed strategy equilibrium that uses strategy 1 with probability $\rho$ with probability $\rho$ and strategy 2 with probability $1-\rho$, the limiting $\phi(u) = cu(\rho-u)(1-u)$ with $c>0$. Here, as in the case $q<1$, the fixed point $\rho$ is attracting.
To translate  Theorem 4 in \cite{CD} to our situation, we note that $w=\ep_L^2$ and $n=L^d$.

\begin{theorem} \label{expsurv} 
Suppose that $\ep_n \sim C n^{-2\alpha/3}$, where $0<\alpha<1$, and that we start from a product measure in which 
each type has positive density. Let $N_1(t)$ be the number of sites occupied by 1's at time $t$. 
There is a $c > 0$ so that for any $\eta > 0$ if $n$ is large and
$\log n \le t \le \exp( c n^{(1-\alpha)})$, then $N_1(t)/N \in (\rho-\eta,\rho+\eta)$ with high probability. 
\end{theorem}

\noindent
The intuition behind the answer is that after space is rescaled the volume of the torus is asymptotically $n^{(1-\alpha)}$. Theorem \ref{expsurv} is a lower bound so it does not rule out survival for time $\exp(cn)$. However, Cox and Durrett proved for the contact process with fast voting introduced by Durrett, Liggett, and Zhang \cite{DLZ} 

\begin{theorem} \label{CPdeath}
There is a $C<\infty$ so that extinction in the contact process plus fast voting occurs by time $\exp(cn^{1-2\alpha/d}\log n) $ in $d \ge 3$.
\end{theorem}

\noindent
Theorem \ref{expsurv} can be generalized to the $q$-voter with $q<1$ since it only relies on the hydrodynamic limit in Theorem \ref{hydro} and a block construction. Theorem \ref{CPdeath} does not extend, because $\xi\equiv 1$ is an absorbing state, and this limits our ability to suddenly kill the process.

\clearp 

\section{Graphical representation, duality} \label{sec:grep}

\subsection{Voter model} \label{ss:vmd}

We begin by describing the graphical representation and duality for the voter model in which the neighbors of $x$ are $x+{\cal N}$ and ${\cal N} = \{ y_1, \ldots y_k\}$. The state of the voter model at time $t$ is $\xi_t : \ZZ^d \to \{0,1\}$ where $\xi_t(x)$ gives the opinion of the individual at $x$ at time $t$. We write $y \sim x$ to indicate that $y$ is a neighbor of $x$.
In the usual voter model, the rate at which the voter at $x$ changes its opinion from $i$ to $j$ is
$$
c^v_{i,j}(x,\xi) = 1_{(\xi(x)=i)} f_j(x,\xi), 
$$
where $f_j(x,\xi) = (1/k) \sum_{i=1}^k 1(\xi(x+y_i)=j)$ is the fraction of neighbors in state $j$.

To study the voter model, it is convenient to construct the process on a {\it graphical representation}, introduced by Harris \cite{H76} and further developed by Griffeath \cite{G78}. For each $x \in \ZZ^d$ and $y\in x+{\cal N}$ let $T_m^{x,y}$, $m \ge 1$, be the arrival times of a Poisson process with rate $1/k$. At the times $T^{x,y}_n$, $n \ge 1$, the voter at $x$ decides to change its opinion to match the one at $y$.  To indicate this, at time $T^{x,y}_n$ we write a $\delta$ at $x$ and draw an arrow from $y$ to $x$. To calculate the state of the voter model on a finite set, we start at the bottom and work our way up. We think of the 1's in the initial configuration as sources of fluid, the $\delta$'s as dams that block the fluid, while the arrows move the fluid in the direction indicated. Arrows from $y$ to $x$ arrive just after the $\delta$. A nice feature of this approach is that it simultaneously constructs the 
process for all initial conditions so that if $\xi_0(x) \le \xi'_0(x)$ for all $x$, then for all $t>0$ 
we have $\xi_t(x) \le \xi'_t(x)$ for all $x$.

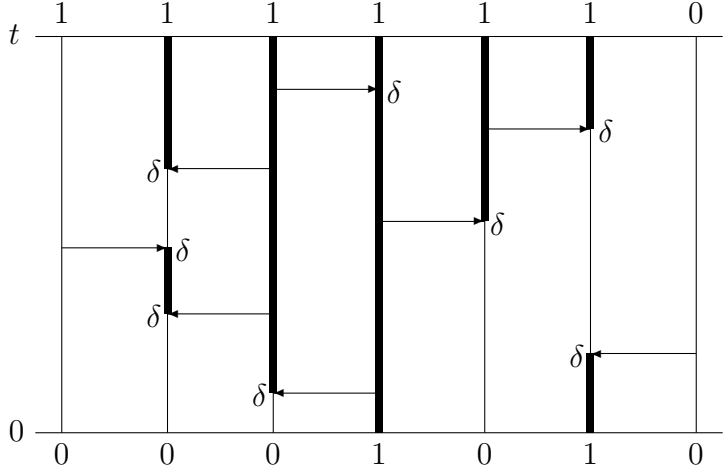
\begin{figure}[ht]
\begin{center}
\begin{picture}(320,220)
\put(30,30){\line(1,0){260}}
\put(30,180){\line(1,0){260}}
\put(40,30){\line(0,1){150}}
\put(80,30){\line(0,1){150}}
\put(120,30){\line(0,1){150}}
\put(160,30){\line(0,1){150}}
\put(200,30){\line(0,1){150}}
\put(240,30){\line(0,1){150}}
\put(280,30){\line(0,1){150}}
\put(37,18){0}
\put(77,18){0}
\put(117,18){0}
\put(157,18){1}
\put(197,18){0}
\put(237,18){1}
\put(277,18){0}
\put(37,185){1}
\put(77,185){1}
\put(117,185){1}
\put(157,185){1}
\put(197,185){1}
\put(237,185){1}
\put(277,185){0}
\put(20,27){0}
\put(20,177){$t$}
\put(120,160){\vector(1,0){40}}
\put(163,155){$\delta$}
\put(200,145){\vector(1,0){40}}
\put(243,140){$\delta$}
\put(120,130){\vector(-1,0){40}}
\put(72,125){$\delta$}
\put(160,110){\vector(1,0){40}}
\put(202,105){$\delta$}
\put(40,100){\vector(1,0){40}}
\put(83,95){$\delta$}
\put(120,75){\vector(-1,0){40}}
\put(72,70){$\delta$}
\put(280,60){\vector(-1,0){40}}
\put(232,55){$\delta$}
\put(160,45){\vector(-1,0){40}}
\put(112,40){$\delta$}
\linethickness{1.0mm}
\put(160,30){\line(0,1){150}}
\put(80,75){\line(0,1){25}}
\put(80,130){\line(0,1){50}}
\put(120,45){\line(0,1){135}}
\put(200,110){\line(0,1){70}}
\put(240,30){\line(0,1){30}}
\put(240,145){\line(0,1){35}}
\end{picture}
\caption{Voter model graphical representation}
\end{center}
\end{figure}

To define the {\it dual process} starting from $x$ at time $t$, we set $\zeta^{x,t}_0 = x$ and work down the graphical representation. A particle stays at its current location until the first time that it encounters a $\delta$. At this point it jumps across the edge in the direction opposite its orientation.  A little thought reveals that the path of a single particle in $\zeta^{x,t}_s$, $0 \le s \le t$, is a random walk that at rate 1 jumps to a randomly chosen neighbor. Intuitively, $\zeta^{x,t}_s$ gives the source at time $t-s$ of the opinion at $x$ at time $t$. That is,
$$
\xi_t(x) = \xi_{t-s}(\zeta^{x,t}_s).
$$
The example in Figure 3 should help explain the definitions. Here we work backwards to determine the states of the two sites marked by `?'. The dark lines indicate the locations of the two dual particles. The family of particles $\zeta^{x,t}_s$ are coalescing random walks. That is, if a particle $\zeta^{x,t}_s$ lands on the site occupied by $\zeta^{y,t}_s$, the two particles coalesce to form a single particle, and we know that $\xi_t(x)=\xi_t(y)$.

\begin{figure}[ht]
\begin{center}
\begin{picture}(320,220)
\put(30,30){\line(1,0){260}}
\put(30,180){\line(1,0){260}}
\put(40,30){\line(0,1){150}}
\put(80,30){\line(0,1){150}}
\put(120,30){\line(0,1){150}}
\put(160,30){\line(0,1){150}}
\put(200,30){\line(0,1){150}}
\put(240,30){\line(0,1){150}}
\put(280,30){\line(0,1){150}}
\put(37,18){0}
\put(77,18){0}
\put(117,18){0}
\put(157,18){1}
\put(197,18){0}
\put(237,18){1}
\put(277,18){0}
\put(77,185){?}
\put(197,185){?}
\put(20,27){0}
\put(20,177){$t$}
\put(120,160){\vector(1,0){40}}
\put(163,155){$\delta$}
\put(200,145){\vector(1,0){40}}
\put(243,140){$\delta$}
\put(120,130){\vector(-1,0){40}}
\put(72,125){$\delta$}
\put(160,110){\vector(1,0){40}}
\put(202,105){$\delta$}
\put(40,100){\vector(1,0){40}}
\put(83,95){$\delta$}
\put(120,75){\vector(-1,0){40}}
\put(72,70){$\delta$}
\put(280,60){\vector(-1,0){40}}
\put(232,55){$\delta$}
\put(160,45){\vector(-1,0){40}}
\put(112,40){$\delta$}
\linethickness{1.0mm}
\put(80,180){\line(0,-1){50}}
\put(120,130){\line(0,-1){85}}
\put(160,110){\line(0,-1){80}}
\put(200,180){\line(0,-1){70}}
\end{picture}
\caption{Dual coalescing random walk}
\end{center}
\end{figure}
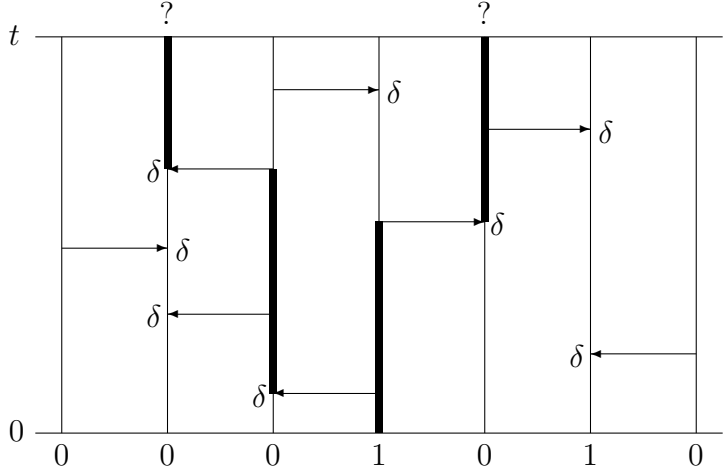

To illustrate the power of duality, we analyze the asymptotic behavior of the voter model on $\ZZ^d$, proving a result of Holley and Liggett \cite{HL}. 
In dimensions 1 and 2, nearest neighbor random walk is recurrent, so the voter model clusters, i.e.,
$$
P( \xi_t(x) \neq \xi_t(y) ) \le P( \zeta^{x,t}_t \neq \zeta^{x,t}_t ) \to 0.
$$  
In $d \ge 3$ random walks are transient so differences in opinion persist as $t \to \infty$. 
Let $\xi^u_t$ be the voter model starting from product measure in which 1's have density $u$, i.e., the initial voter opinions are independent and $=1$ with probability $u$. For a finite set $B \subset \ZZ^d$, let $\zeta^{B,t}_s = \cup_{x\in B} \zeta^{x,t}_s$. The distribution of  $\zeta^{B,t}_s$ does not depend on $t$ so we drop the superscript $t$. Duality implies
$$
P( \xi^u_t(x) \equiv 0 \hbox{ on } B ) = P( \xi^u_0(y) = 0 \hbox{ for all } x \in \zeta^{B}_t ) = E\left( (1-u)^{|\zeta^{B}_t |} \right)
$$
As $t \uparrow \infty$ , $|\zeta^{B}_t | \downarrow |\zeta^{B}_\infty |$. From this it follows that
\beq
P( \xi^u_t(x) \equiv 0 \hbox{ on } B ) \to  E\left( (1-u)^{|\zeta^{B}_\infty |} \right)
\label{Econv}
\eeq
The probabilities on the left-hand side of \eqref{Econv} are enough to determine the distribution of the limit $\xi^u_\infty$. Since the limit exists, it is a stationary distribution that we denote by $\nu_u$.

Before moving on, we note that the duality equation can be written as
\beq
P( \xi^A_t \cap B \neq \emptyset ) = P( A \cap \zeta^B_t \neq \emptyset )
\label{adual}
\eeq
where $\xi^A_t$ is the voter model starting with 1's on $A$ and $\zeta^B_t$ is the coalescing random walk starting with particles on $B$. This holds because the left-hand side is the probability of a path from $A \times \{0\}$ up to $B \times \{t\}$, while the right-hand side is the probability of a path from $B \times \{t\}$ down to $A \times \{0\}$.
  There are several types of duality. This one is called {\it additive} because $\xi^{A \cup B}_t = \xi^A_t \cup \xi^B_t$, a property that holds because $\xi^A_t$ is defined to be the set of sites at time $t$ that can be reached from a path starting in $A$.

\subsection{Nonlinear voter models} \label{ss:nlv}

Though it is tempting to try to find a duality like the one between the voter model and coalescing random walk to help analyze the $q$-voter model, in this section we will prove

\mn
{\bf Claim.} {\it Using the graphical representation described in the previous section we cannot construct a voter model in which the
flip rates depend only on the number of neighbors with the opposite opinion $n_x$ and are nonlinear.} 

\begin{proof}
For simplicity, we only prove the result when the neighborhood has size 4. Consulting Griffeath's book we see that the only gadgets than can be used in the graphical representation are combination of arrows and $\delta$'s. To begin, we will consider the set of processes that can be constructed by only using gadgets that have a $\delta$ at $x$ and a number of arrows that point to x from its neighbors. We call these objects arrow-$\delta$s. Since the flip rates only depend on the number of sites, all arrow-$\delta$s with $k$ arrows have the same rate, $a_k$.

\begin{itemize}
  \item 
When there is a 1 at $x$ the $\delta$ will cause the 1 to flip to a 0. However, the site will only stay a 0 if \emph{all} neighbors connected to $x$ by arrows are in state 0.
\item
When there is a 0 at $x$ then the $\delta$ does nothing, and the site will flip to 1 if there is \emph{at least one} neighbor in state 1 connected to $x$ by an arrow. 
\end{itemize}

The number of $k$-arrow gadgets is $\binom{k}{2}$ so  the flip rates are as follows 

\begin{center}
\begin{tabular}{ccc}
$n_x$ & rate $1\to 0$ & rate $0 \to 1$ \\
0 & 0 & 0 \\
1 & $a_1$ & $a_1 + 3a_2 + 3a_3 + a_4$ \\
2 & $2a_1 + a_2$ & $2a_1 + 5a_2 + 4a_3 + a_4$ \\
3 & $3a_1 + 3a_2 + a_3$ & $3a_1 + 6a_2 + 4a_3 + a_4$ \\
4 & $4a_1 + 6a_2+4a_3+a_4$ & $4a_1 + 6a_2+4a_3+a_4$ 
\end{tabular}
\end{center}   
\noindent

\noindent
If we add $\delta$'s with no arrows then they will flip 1s even when all their neighbors are 1. 
If $a_2$, $a_3$, or $a_4$ is positive the rate of flipping $1 \to 0$ is $<$  the rate of flipping $0 \to 1$.
when $n_x = 1,2,3$. Adding arrows with no $\delta$s will only further increase the rates of flips $0 \to 1$. 
\end{proof}

\subsection{Duality for voter model perturbations} \label{ss:vpdual}

In the previous section we have shown that the $q$-voter does not have an additive dual. In this section we will introduce a generalization of the graphical representation used in Section \ref{ss:vmd} that allows us to construct voter model perturbations. This idea goes back to \cite{DN}. See also Section 2 in \cite{CDP}. Calculating the state of the process is not as simple as in the additive case, but it does allow us to compute the state of the process on a finite set $B$ at time $t$ by working backwards from time $t$.

Voter model perturbations have flip rates 
\beq
c^\delta_{i,j} = f_j + \delta^2 h_{i,j}(x,\xi),
\label{frates}
\eeq
where $f_j$ is the fraction of neighbors in state $j$.
The perturbation function $h_{ij}$, $j\ne i$, may be negative (and this happens when $q>1$) but 
in order for the analysis in \cite{CDP} to work, there must be a 
law $q$ of $(Y_1, \ldots Y_k) \in (\ZZ^d)^k$ and a functions $g_{i,j} \ge 0$, so that for some $\gamma < \infty$, we have 
\beq
h_{i,j}(x,\xi) = - \gamma f_j + E_{Y}[g_{i,j}(\xi(x+Y_1), \ldots \xi(x+Y_k))].
\label{vptech}
\eeq
In our situation $Y_1, \ldots Y_k$ are $k$ neighbors in ${\cal N}$ and $g_{i,j}$, which does not depend on $\ep$, is the fraction of sites $x+Y_1, \ldots x+Y_k$ in state $j=1-i$ raised to the $q$th power.

Suppose now that we have a voter model perturbation of the form \eqref{frates} which satisfies \eqref{vptech}. 
We construct the voter model portion as in Section \ref{ss:vmd}.
We call the arrow-$\delta$s {\bf voter events}. To add the perturbation we let 
$$
\|g_{i,j}\| = \sup_{\eta \in \{0,1\}^M} g_{i,j}(\eta_1, \ldots \eta_k)
$$ 
and introduce Poisson processes
$T^{x,i,j}_m$, $m\ge 1$ with rate $r_{i,j} = \ep \|g^\ep_{i,j}\|$, where $\ep=\delta^2$, and independent random variables $U^{x,i,j}_m$, $m\ge 1$ uniform on $(0,1)$. At the times $t=T^{x,i,j}_m$ with $m \ge 1$ we draw arrows from $x+Y^i$ to $x$ for $1\le i \le k$. 
We call this a {\bf branching event.} If $\xi_{t-}(x)=i$ and 
\beq
r_{i,j}U^{x,i,j}_k < g_{i,j}(\xi_{t-}(x+Y_1), \ldots \xi_{t-}(x+Y_k))  
\label{jrule}
\eeq
then we set $\xi_t(x)=j$. The uniform random variables slow down the transition rate from the maximum possible rate $r_{i,j}$ to the one appropriate for the current configuration.

To define the dual, we proceed as before. When a particle encounters a $\delta$ associated with a voter event, it jumps to the other end of the arrow. When a particle encounters the head of an arrow associated with a branching event it gives birth to new particles at the other ends of all of the arrows. If either action results in two particles on the same site they coalesce to 1. Let $I^{B,t}_s$ be the set of particles at time $t-s$ when we start with particles on $B$ at time $t$. Durrett and Neuhauser \cite{DurNeu} called $I^{B,t}_s$ the {\bf influence set} because 

\begin{lemma}
If we know the values of $\xi_{t-s}$ on $I^{B,t}_s$, then using the graphical representation (including the associated uniform random variables) we can compute the values of $\xi_t$ in $B$ by working our way up the graphical representation starting from time $t-s$ and determining the changes that should be made in the configuration at each jump time.
\end{lemma} 

This fact should be clear from the construction. A formal proof can be found in Section 2.6 of \cite{CDP}.
The {\bf computation process}, as it is called in \cite{CDP}, is complicated, but is useful because
up to time $t/\ep_n$ there will only be $O(1)$ branching events affecting particles in the dual. 
 
\clearp 

\section{Prolonged persistence} \label{sec:proper}

In this section, we will prove Theorem \ref{persist}. The key is to bound the difference between the density of the particle system and the ODE, using a result of Darling and Norris \cite{DN}. Section \ref{ss:DarNor} describes this result and the work needed to apply it to finish the proof of Theorem \ref{persist}. Sections \ref{ss:igbr}, \ref{ss:dens}, \ref{ss:TvsT}, and \ref{ss:diffm} complete this work and Section \ref{ss:fd} gives the final details.

\subsection{Darling-Norris theorem} \label{ss:DarNor}

To state the result from \cite{DN} result we need to introduce some notation.  Let $\xi_t$ be a continuous time Markov chain with countable state space $S$ and jump rates $q(\xi,\xi')$. In our case $\xi_t$ will be the state of the $q$-voter model on the torus. We are interested in proving an ODE limit for $X_t =x(\xi_{t/\ep_n})$ where
$$
x(\xi_{t/\ep_n}) = \frac{1}{n} \sum_{x\in \TT_n} \xi_{t/\ep_ n}(x).
$$
For each $\xi\in S$ we define the infinitesimal drift 
$$
\beta(\xi) = \sum_{\xi'\neq\xi} (x(\xi')-x(\xi)) q(\xi,\xi')
$$
We let $b$ be the drift of the proposed deterministic limit $x_t$. In our case
$$
x_t = x_0 \pm \int_0^t b(x_s) \,ds, \,\,\,\,\,\,\, \quad\hbox{$b(x) = c_k x(1-x)(1-2x)f_k(x)$},
$$
where $f_k(x)$ is a polynomial with $f_k(0)=f_k(1)=1$ that is positive on $[0,1]$ and only depends on the number of neighbors $k$ . The sign is $+$ for $q=1-\ep_n$ and $-$ for $q=1+\ep_n$. The crucial theorem from \cite{DN} is

\begin{theorem} \label{DarNor}
For each fixed $t_0$ and $\eta  > 0$,
$$
P \left( \sup_{s \le t_0} |X_s - x_s| > \eta \right) \le 2e^{-\gamma^2/(2At_0)} + P( \Omega_0^c \cup  \Omega_1^c \cup  \Omega_2^c )
$$ 
\end{theorem}

To make this statement meaningful we need more definitions. To measure the size of the jumps we let
$\sigma_\theta(y) = e^{\theta|y|} - 1 - \theta|y|$ and let
$$
\phi(\xi,\theta) = \sum_{\xi'\neq \xi} \sigma_\theta( x(\xi')-x(\xi) )q(\xi,\xi').
$$
The good sets $\Omega_i$, $i=0,1,2$ are given by
\begin{align}
&\Omega_0 = \{ |X_0 - x_0| \le \gamma \} \label{om1} \\
&\Omega_1  = \left\{ \int_0^{t} |\beta(\xi_{s/ \epsilon_n}) - b(X_s)| \, ds \le \gamma \right\},\label{om2} \\
&\Omega_2   = \left\{ \int_0^t \phi(\xi_{s/ \epsilon_n}, \theta) \, ds \le \theta^2At/2 \right\}. \label{om3}
\end{align}
The parameters in these events are coupled by the following relationships.  If we let $K$ be the Lipschitz constant
of the drift $b$ and $\eta$ be the upper bound on the error in the approximation by the differential equation in Theorem \ref{DarNor},
then 
$$
\gamma = \eta e^{-Kt_0}/3 \quad\hbox{and}\quad \theta = \gamma/(At_0), \quad\hbox{where $A>0$}. 
$$

It is clear that our $b(x)$ is Lipschitz continuous.
Our assumption that $U_n(0) \to u_0$ implies that $\Omega_0^c = \emptyset$ for large $n$. 
To bound $P(\Omega_2^c)$, we will choose an $A > 0$ that works well. We begin with a useful lemma:

\begin{lemma} \label{PoiLD}
If $Z \sim \hbox{Poisson}(\lambda)$, then 
$$
P(Z \ge 2\lambda) \le \exp(-\gamma(2)\lambda)
$$
where $\gamma(2)$ is a constant independent of $\lambda$.
\end{lemma} 

\begin{proof} The moment generating function of $Z$ is
$$
E\exp(\theta Z) \le \exp(\lambda (e^\theta-1 )).
$$
Taking $\theta = \log 2$, we have $E\exp(Z \log 2) = \exp(\lambda)$, so using Chebyshev's inequality we have
$$
P( Z \ge 2\lambda ) \le \exp(- 2 \log 2 \lambda), \exp(\lambda)  
$$
which proves the result with $\gamma(2)=2\ln 2 - 1$. 
\end{proof}

The process $X_t$ has jumps of size $1/n$ at total rate $n/\ep_n$. As $\theta |y| \to 0$, we have $\sigma_\theta(y) \sim \theta^2y^2/2$. So, when $\theta |y|$ is small, $\sigma_\theta(y) \sim \theta^2y^2$. Using Lemma \ref{PoiLD}, the probability of $2t_0n/\ep_n$ jumps during time $[0,t_0]$ is $\le \exp(-\gamma(2)t_0 n/\ep_n)$. When this occurs, and $n$ is large, the integral in $\Omega_2$ is 
$$
\le \frac{\theta^2}{ n^2} \cdot \frac{2 t_0 n}{\ep_n} = \theta^2 t_0 \cdot \frac{2}{n \ep_n}.
$$ 
Thus, for the event $\Omega_2$ to hold, we need $2/(n \ep_n) \ll A/2$.  Since $\ep_n \sim Cn^{-a}$ with $2/3< a < 1$, we have

\begin{lemma} If $t_0$ and $\gamma$ are fixed and $A=n^{-(1-a)/3}$ then $e^{-\gamma^2/(2At_0)} \to 0$ and $P(\Omega_2^c) \to 0$ exponentially
fast as $n\to\infty$.   
\end{lemma}

\subsection{Ignoring branching} \label{ss:igbr}

The remainder of Section \ref{sec:proper} is devoted to bounding $P(\Omega_1^c)$.
To begin to do this, we return to the original time scale. We define $\tilde{\xi}_s$ to be the same as $\xi_s$ at time $s = t/\ep_n- n^{b}$, while on the time interval  $[t/\ep_n - n^b, t/\ep_n ]$, $\tilde{\xi}_s$ only has voter events, ignoring the perturbation. The value $b \in (2/3,a)$ is chosen so that lineages in the dual coalescing random walk will have time to wrap around the torus but, as we will now show, the perturbation will not have much effect.  Let 
$$
\tilde{X}_t= \frac{1}{n} \sum_{x \in \TT_n} \tilde{\xi}_{t/\ep_n}(x)
$$ 
be the density of this new process $\tilde{\xi}$. 

We will now show that ignoring the perturbation changes the values of more that $\eta n$ sites with a stretched exponentially small probability. 

\mn
{\bf Step 1.} The number of perturbation events $M$  in time $n^b$ is bounded by a Poisson($\lambda$) random variable with $\lambda= Cn^{1+b+a}$. 
Lemma \ref{PoiLD} implies that 
\beq
P(M \ge 2 \lambda ) \le \exp(-\gamma(2)\lambda) \le \exp(-C\gamma(2)n^b),
\label{bdpert}
\eeq
since $\lambda \ge Cn^b$.

\mn
{\bf Step 2.} Let $\eta_t(x) = |\xi_t(x)-\tilde\xi_t(x)|$, so that $\eta_t(x)=1$ means there is a discrepancy between the two processes $\xi_t$ and $\tilde\xi_t$ at position $x$. We want to prove that $\sum_x \eta_{t/\ep_n}(x)$ is less than $\eta n$  with a stretched exponentially small probability. To do this, note that when an edge $(x,y)$ with $\eta_s(x)=0$ and $\eta_s(y)=1$ is hit by a voter event (that is, there is an arrival in the Poisson process $T^{x,y}$ or $T^{y,x}$), then the 1 is changed to $0$ with probability 1/2 (when the arrival is in $T^{x,y}$) and the 0 is changed to a 1 with probability 1/2 (when the arrival is in $T^{y,x}$). Thus, the change in the number of discrepancies due to voter events is a martingale. The change is always $\le 1$ so if there are $N$ jumps, then by Azuma's inequality
$$
P( |X_N - X_0| \ge z | N = n_0) \le 2 \exp(-z^2/2n_0)
$$
If $N$ is the number of changes due to voter events in the time interval $[t/\ep_n -  n^b, t/\ep_n]$, then $ N \le \hbox{Poisson}(n^{b+1})$. By Lemma \ref{PoiLD},
$$
P(N \ge 2 n^{1+b}) \le \exp(-\gamma(2)n^{1+b}).
$$ 
Note that if $n_0 < 2n^{1 + b}$, then $2 \exp(-z^2/2n_0) < 2 \exp( -z^2/4n^{1+b})$. So, taking $z=\eta n$ and $N=2n^{1+b}$, we get
\begin{align}
P( |X_N - X_0| \ge \eta n) \le 2\exp(-\eta^2 n^{1-b}/4).
\label{bdpert2}
\end{align}

\subsection{Bounding the density} \label{ss:dens}

The results in the previous section show that on the interval $[t/\ep_n - n^b, t/\ep_n]$ we can ignore the perturbation and assume that the process evolves like the voter model. To understand the distribution of 1's at time $t/\ep_n$ we will use results of Bramson and Griffeath \cite{BGrenorm}, and Z\"ahle \cite{Zrenorm}. The first reference only treats $d=3$. The second covers $d \ge 3$ and is more detailed, so we will follow it. 

Let $\zeta^\lambda: \ZZ^d \to \{0,1\}$ have the distribution of the equilibrium of a finite range voter model on $\ZZ^d$ with density $\nu_\lambda$. For an explanation of this and the other basic facts about the voter model that we will use, see Liggett's book \cite{L99}. For simplicity we will do calculations for the nearest neighbor case. The results are the same in the finite range case, but are more awkward to write since, for example, the limiting normal has a general covariance matrix, we cannot use the reflection principle, etc. To formulate the limit theorem in \cite{Zrenorm}, we will write the process at a fixed time as a random field
$$
F_\lambda(\phi) = \sum_{i \in \ZZ^d} [\zeta^\lambda(i) - \lambda] \phi(i),
$$
where $\phi$ is a member of a suitable class of test functions. To rescale space, we let 
$$
F_{\lambda,r}(\phi) = F_{\lambda}(\phi_r) \quad\hbox{where}\quad \phi_r(x) = r^{-(d+2)/2} \phi(x/r).
$$
Theorem 1 on pages 1265--1266 of \cite{Zrenorm} shows that in our nearest neighbor case
$$
F_{\lambda,r}(\phi) \Rightarrow \hbox{Normal}(0, a_d \lambda (1-\lambda) B(\phi,\phi)),
$$
where $\Rightarrow$ denotes weak convergence as $r \to \infty$, $\text{Normal}(\mu,\sigma^2)$ is a one-dimensional normal distribution
with mean $\mu$ and variance $\sigma^2$, and $B$ is the bilinear function
$$
B(\phi,\psi) = \iint \frac{\phi(x) \psi(y)}{|x-y|^{(d-2)/2}} \, dx \, dy.
$$
Restricting our attention now to $d=3$, Z\"ahle's result implies that
\beq
\widehat S_r \equiv [\lambda (1-\lambda)]^{-1/2} r^{-5/2} \sum_{x \in [-r/2,r/2]^3} \left[\zeta^\lambda (x) -\lambda \right] 
\Rightarrow \hbox{Normal}(0, c_{3,\lambda} ) 
\label{boxsum}
\eeq
Bramson and Griffeath \cite{BGrenorm} prove \eqref{boxsum} by the method of moments, which gives \
\beq
E(\widehat S_r)^{2m} \to c_{3,\lambda}^{2m}  \mu_m \quad\hbox{where} \quad \mu_m = (2m-1)(2m-3) \cdots 3 \cdot 1.
\label{bsmom}
\eeq
In our situation, we need a slightly different result. In particular, these results are for the voter model on $\ZZ^3$, and we need a result for the voter model on the 3-d torus. Let
$$
T_r \equiv [\lambda (1-\lambda)]^{-1/2} r^{-5/2} \sum_{x \in Q(r)} [\xi_{t/\ep_n}(x) -\lambda]
$$
where $\lambda$  is the fraction of sites in state 1 at time $t/\ep_n-n^b$, and $Q(r)$ is a fixed cube with side $r=n^\beta$ with $\beta< 1/3$. To prove a limit result for $T_r$ we will sandwich it between $\widehat S_r$ and
$$
\bar S_r \equiv [\lambda (1-\lambda)]^{-1/2} r^{-5/2} \sum_{x \in Q(r)} [\bar\zeta^\lambda_{\sigma(n)} (x) -\lambda ],
$$
where $\hat\zeta^\lambda_{\sigma(n)}$ is the voter model on the torus starting from product measure with density $\lambda$ and run for time $\sigma(n)=n^{0.6}$. To couple this with $T_r$ we create $\bar S_r$ by running coalescing random walks starting at time $t/\ep_n$ from points in $Q(r)$ backwards in time for $\sigma(n)$, and then use independent coin flips with probability $\lambda$ of heads (1) and $1-\lambda$ of tails (0) to determine the states of the sites.

\mn
(i) {\bf With stretched exponentially small probability, no coalescing random walk will move more than $n^{0.33}$ in any coordinate by time $\sigma(n)=n^{0.6}$.} 

\mn
\begin{proof} We will use a special case of (7.3) on page 553 in Feller volume II \cite{Feller}.

\begin{lemma} \label{tailbound}
Let $w_1, w_2, \ldots w_k$ be i.i.d.~with $P(w_i=1)=P(w_i=-1)$. Then if
$W_k = w_1 + \cdots w_k$, $\ep>0$, and $x=o(k)$, we have
$$
P( W_k/\sqrt{k} \ge x) \le \exp(-(1-\ep) x^2/2 ). 
$$
\end{lemma}

\noindent
Taking $k=n^{0.6}$ and $x=0.03$ , it follows that the probability some coalescing random walk starting inside the cube $Q(r)$ and run for time $\sigma(n)$ moves by more than $n^{0.33}$ in  any coordinate is
$$
\le 2 \cdot 6r^3 \exp( -(1-\ep) n^{0.06}/2 ).
$$
Here the 2 comes from using the reflection principle to relate the maximum to the value at time $n^{0.6}$, and 6 is 3 coordinates times 2 signs.
\end{proof}

\medskip
The result (i) implies that with very high probability there is no difference between the coalescing starting from $Q(r)$ with $r = n^\beta$ for $\beta < 1/3$, run to time $\sigma(n) = n^{0.6}$ on the torus or on $\ZZ^3$. 

\mn
(ii) {\bf There is a $\gamma>0$ so that at all times $t \ge (k+1)n^{2/3}$, the total variation between the distribution of a nearest neighbor random walk on the torus and the uniform distribution is  $\le (1-\gamma)^{k}$.}

\begin{proof} 
To prove the result, we use a simple coupling. At time $n^{2/3}$ the distribution of each particle has a density that is $\ge \gamma/n$ at each point of the torus.  At time $n^{2/3}$ the distribution has the form $\gamma \cdot \mu_n + (1-\gamma) q_n$, where $\mu_n$ is uniform on the torus and $q_n$ is some transition probability. Uncoupled mass at time $(k-1)n^{2/3}$ can be coupled to the uniform distribution with probability $\ge \gamma$ at time $kn^{2/3}$ and the desired result follows. 
\end{proof} 

\mn
{\bf Definition of $T'_n$.} We continue the construction of $T_r$: from the end of the construction of $\bar S_r$ at time $\sigma(n)$, we run the coalescing random walk particles on $\ZZ^3$. To assign values to the lineages at time $n^b$ we extend the configuration on the torus at that time to be periodic on $\ZZ^3$. It follows from (ii) that with very high probability there is no difference between flipping coins at time $n^{0.6}$ to determine the states of the sites in the sum $\bar S_n$ or continuing to run the coalescing random walks on $\ZZ^3$ until time $n^b$. Having done this, we no longer perfectly reproduce $T_n$, so we call the result $T'_n$. The good news is that when we run the  coalescing random walk on $\ZZ^3$ starting at $\sigma(n)$, we will have $T'_r \prec \widehat S_n$. That is, the coalescing random walk clusters in $T_r'$ are contained in clusters in $\widehat S_r$.

\mn
{\bf To prove the result in \eqref{boxsum}}, Z\"ahle defines a \emph{cluster} to be a set of sites that coalesce to the same limiting particle, and lets $Z_{r,k}$, $1 \le k \le K(r)$ be the cluster sizes and lets $\eta_{r,k}$ be independently $=1$ with probability $\lambda$ and $=0$ with probability $1 - \lambda$. As she notes in (3.6) on page 1274,
\beq
\widehat S_r =_d r^{-5/2} \sum_{k=1}^{K(r)} Z_{r,k} \cdot (\eta_{n,k} - \lambda)
\label{clrep}
\eeq
If we condition on the $Z_{r,k}$, then we have a sum of independent random variables. If we let $v_n^2 = \sum_k Z_{n,k}^2$, then using Lyapunov's theorem (see the bottom of page 1275) it follows that
$$
(\hat S_r/v_n \mid {\cal Z}) \Rightarrow \chi,
$$
where ${\cal Z}$ is the $\sigma$-field generated by the $Z_{r,k}$ and $\chi$ is a standard normal. In Lemma 1 on page 1276 in \cite{Zrenorm} she shows that $v_n^2$ converges in probability to a constant, so if we remove the conditioning we get the same limit. Lemma 2 computes the limit of $Ev_n^2$ and \eqref{boxsum} follows. 

\mn
{\bf The last argument can be applied to $\bar S_n$} to conclude that it converges to a normal distribution. To find the limiting variance we compute
$$
\sum_{x,y\in Q(r)} E( \bar\zeta^\lambda_{\sigma(n)} (x) -\lambda )( \bar\zeta^\lambda_{\sigma(n)} (y) -\lambda )
$$
When the coalescing random walks starting from $x$ and $y$ do not coalesce, the states at $x$ and $y$ are independent; otherwise, they are equal. Thus, if we let $\tau_{x,y}$ be the time the two coalescing random walks hit, then the above sum is
$$
\sum_{x,y\in Q(r)} \lambda(1-\lambda) P( \tau_{x,y} \le n^{0.6} )
$$
Using the local central limit theorem,
$$
P( n^{0.6} \le \tau_{x,y} < \infty ) \approx 2\beta_d \int_{n^{0.6}}^\infty \frac{1}{(2\pi t)^{3/2}}\exp \left( - |x-y|^2/2t \right) \, dt
$$
The right-hand side gives the expected amount of time the two particles spend together. When they hit they spend an exponential rate 2 amount of time together. In addition, they will hit a geometric number of times with success probability $\beta_d$.  
Changing variables $t = |x-y|^2/2s$, $dt = -|x-y|^2/(2s^2)$ the integral becomes
\begin{align*}
& \int_0^{|x-y|^2/n^{0.6}} \left( \frac{s}{\pi|x-y|^2} \right)^{3/2} e^{-s} \left( \frac{\pi|x-y|^2}{2s^2} \right) \, ds \\
&=  \frac{1}{2\pi^{3/2} |x-y|} \int_0^{|x-y|^2/n^{0.6}} s^{-1/2} e^{-s} \, ds  \le C n^{-0.3}.
\end{align*}
Consulting Lemma 4 in \cite{Zrenorm}  we find
$$
P( \tau_{x,y} < \infty ) \sim c'_3/|x-y|
$$
Using the formula for $c'_3$ it follows that the asymptotic variance for $\bar S_r$ is the same as for $\widehat S_r$.

\mn
{\bf Limit theorem for $T'_r$.} Let $X_{r,k} \prec Y'_{r,k} \prec Z_{r,k}$ be the cluster sizes in $\bar S_r$, $T'_r$, and $\widehat S_r$. The limiting variances of the unnormalized sums are
$$
\sum_k EX_{r,k}^2 \le \sum_k E(Y'_{r,k})^2 \le \sum_k EZ_{r,k}^2 
$$
Since the top and bottom sums have the same asymptotics, this gives us the Gaussian limit theorem for $T'_n$. Replacing 2 by $2m$ and recalling that Bramson and Griffeath \cite{BGrenorm} proved their result for $\widehat S_r$ by the method of moments gives the desired results for $T'_r$:
\begin{align}
&T'_r \equiv [\lambda (1-\lambda)]^{-1/2} r^{-5/2} \sum_{k=1}^{K(r)}Y'_{r,k} (\eta_{r,k} - \lambda)  
\Rightarrow {\cal N}(0, c_{3,\lambda} ) 
\label{boxsumT}
\\
& E(T'_r)^{2m} \to c_{3,\lambda}^{2m} (2m-1)(2m-3) \cdots 3 \cdot 1
\label{bdmomT}
\end{align}

The last result implies
\beq
r^{2m\beta} P( |T'_r| \ge r^\beta ) \le E(\bar S_r)^{2m} \to E\chi^{2m}
\label{Chebym}
\eeq
so if we let $\tilde D'_r = [\lambda (1-\lambda)]^{1/2} r^{5/2} T'_r$, (i.e., we remove the scaling) then  
\beq
P( |D'_r| \ge [\lambda (1-\lambda)]^{1/2} r^{5/2+\beta} )  \le C_m r^{-2m\beta}.
\label{mombdT}
\eeq
This is the concentration result we desired for $T'_n$. Recall that $T'_n$ was constructed as a slight modification of $T_n$, which is the true rescaled and centered density that we which to prove results about.

\subsection{Controlling the difference between $T'_n$ and $T_n$} \label{ss:TvsT}

The goal in this section is to generalize \eqref{mombdT} to $T_r$.

\medskip
{\bf Bounding the number of extra coalescences in $T'_n$.}
When we went from the torus to $\ZZ^3$ we may have eliminated some coalescence in $T_n$ at times in $[n^{0.6}, n^b]$. For this to happen the difference in two particles positions must have wrapped around the torus, an event we call $G$, and the particles projected back to the torus must have hit, an event we call $H$. To bound this event we note that 
$$
P(G \cap H) \le \min\{P(G),P(H)\}.
$$ 
Let $\alpha = 2(1-\ep)/3$. Lemma \ref{tailbound} implies that the probability $G$ happens during $[n^{0.6},n^\alpha]$ is $\le \exp(-n^\eta)$ for some $\eta>0$. On $[n^{\alpha},n^b]$, the probability that a random walk is at a fixed site is $\le 1/n^{1-\ep}$. Thus, for a fixed pair of particles,
$$
P(H) \le C n^b /n^{1-\ep}.
$$

If $r=n^{b(2)/3}$, then $n^{b(2)}$ is a trivial upper bound for the number of particles at time $\sigma(n)$, which holds with probability 1. We will now estimate the number of collisions of a fixed particle with all of the others. This number is increased if we ignore coalescence, and run the particles as independent. We do this so that

\begin{lemma} \label{mpairs}
If $m \ge 1$ and a particle belongs to a cluster of size $2m$ or $2m+1$ with $m \ge 1$ formed by coalescence during $[n^\alpha,n^b]$, then there are at least $m$ disjoint pairs of particles that have coalesced.
\end{lemma}

\begin{proof} Recall that on this time interval we are running the lineages on $\ZZ^3$. We will prove the result by induction. To be able to disentangle the graph constructed by coalescence we will number the particles. Once two particles hit the two future trajectories could be assigned to either particle so we allow ourselves the liberty of be exchanging the labels at any collision. If the cluster has size 2 or 3, this is trivial. Suppose now that $m \ge 2$. Locate the time $t_0$ at which the first two particles coalesced. Call them $x$ and $y$ and let $t_1$ be the first time after $t_0$ that the coalesced particle collided with another one that we call $z$. Remove the $Y$-shaped part of the genealogy leading from $x$ and $y$ to the coalescence at time $t_1$.  Label the lineage coming out $t_1$ the same as the one coming in on $z$'s trajecctory. We have identified one pair of coalescing particles and reduced the number of sites in the cluster by 2, so the result follows by induction.  
\end{proof}     

Given Lemma \ref{mpairs}, our next task is to estimate the probability that $m$ disjoint pairs will coalesce. Using the trivial upper bound $n^{b(2)}$ on the number of lineages, the number of coalescing pairs is
$$
N \le \hbox{Binomial}(n^{2 b(2)}, Cn^b/n^{1-\ep}).
$$
Note that this bounds the number of coalescing pairs that coalesce in the system, not just those that form one cluster.
The expected number is $Cn^{b + 2b(2) + \ep -1}$, where $b$ is larger than 2/3 and can be assumed to be $\le 0.7$. If $b(2)\le 0.1$, then $-\nu = b + 2b(2) + \ep -1 < 0$ when $\ep<0.5$. In this case,
$$
P(N=k) \le \binom{ n^{2b(2)} }{k} (  Cn^{b+\ep -1} )^k \le \frac{C^kn^{-k\nu}}{k!},
$$
so summing gives
\beq
P(N \ge k ) \le e^C n^{-k\nu}.
\label{Kbdd}
\eeq

\mn
{\bf Bounding the size of clusters in $\hat S_r$.}
Formula \eqref{bsmom} tells us that  
$$
E(\widehat S_r)^{2m} \to c_{3,\lambda}^{2m} \mu_m
$$
Using \eqref{clrep} we have
$[(1-\lambda)\lambda^{2m} + \lambda(1-\lambda)^{2m}] \sum_{k=1}^{K(r)} Z_{r,k}^{2m} \le E(\widehat S_r)^{2m}$.
From this we see that when $r$ is large
$$
r^{11m/2} P( \max_k Z_{r,k} \ge r^{5.5/2} ) \le C_m r^{5m}
$$
so we have
\beq
 P( \max_k Z_{r,k} \ge r^{5.5/2} ) \le C_{m,\lambda} r^{-m/2}.
\label{maxZbd}
\eeq
Combining \eqref{Kbdd} and \eqref{maxZbd} we see that if $Y_{r,k}$ are cluster sizes in $T_n$, then
\beq
 P\left(  \max_k Y_{r,k} \ge \frac{m}{2\nu} r^{5.5/2} \right) \le C_{m,\lambda} r^{-m/2}
\label{maxYbd}
\eeq

Combining \eqref{Kbdd} with $k=m/2\nu$ and \eqref{maxYbd} we see that the combined size of the clusters in $T'_n$ but not in $T_n$ is
\beq
\le \frac{m}{2\nu} n^{5.5/2}  \quad\hbox{with probability $1 - C_{m,\lambda} r^{-m/2}$.}
\label{Tdiffbd}
\eeq 
Using this with \eqref{mombdT} and letting $D_r = [\lambda(1-\lambda]^{1/2} r^{5/2} T_n$ it follows that
\beq
P( |D_r| \ge  [\lambda(1-\lambda)]^{1/2} r^{5/2+\beta} \le C_{m,\lambda} r^{-m\beta/2}.
\label{mombdT2}
\eeq

Suppose $r=n^{b(2)/3}$ where $0< b(2)<1$, then
$$
P\left( |D_r| \ge [\lambda(1-\lambda)]^{1/2} n^{5b(2)/6+\beta} \right)  \le C_m n^{-m\beta b(2)/6}
$$
Now, partition the torus into cubes of side $n^{b(2)/3}$. Letting $N_{i}$ be the number of 1's in the $i$th cube we have
$$
P\left( \left| \frac{N_{r,i}}{n^{b(2)}} - \lambda \right| 
\ge [\lambda(1-\lambda)]^{1/2} n^{-b(2)/6+\beta} \right)  \le C_m n^{-m\beta b(2)/6}.
$$
For fixed $\beta>0$, given a $k < \infty$ we can pick $m$ large enough then the right hand side is $\le n^{-(1-b(2))-k}$. 
Then we have,  
\beq 
P\left( \hbox{ for some $i$ }\left| \frac{N_{i}}{n^{b(2)}} - \lambda \right| \ge [\lambda(1-\lambda)]^{1/2} n^{-b(2)/6+\beta} \right)  \le n^{-k}.
\label{mombd3}
\eeq

\subsection{Bounding the difference in the drifts} \label{ss:diffm}

Thus far we have been concerned with the overall density of particles on the torus. However, to successfully bound $P(\Omega^c_1)$ we need to show that if $u$ is the density of ones in the voter model at time $t/\epsilon_n - n^b$, then the empirical finite dimensional distributions on the torus  are close to those of the voter model equilibrium $\nu_u$ at time $t/\epsilon_n +s_n$, where
\beq
s_n = n^{(2+\beta)b(2)/3}.
\label{sndef}
\eeq
The reasoning for introducing this extra time $s_n$ is described below. For $x, y_1, \ldots y_k \in \ZZ^d$ and $v_0, v_1, \ldots v_k \in \{0,1\}$ fixed we let
$$
G_{x,y,v} = \{ \xi(x)=v_0, \xi(x+y_1) = v_1, \ldots \xi(x+y_k) = v_k \}
$$
be a finite dimensional event. For simplicity, we do not display the dependence on the sites $y$ and the states $i$. 

The first step is to partition the torus at time $t/\ep_n$ into boxes with side $r=n^{b(2)/3}$. Using \eqref{mombd3}, we can conclude that with high probability the density in each box is close to $u$, the density of 1's at time $t/\ep_n - n^b$. 
We divide the torus at time $t/\ep_n + s_n$ into cubes with side $n^{b(3)/3}$, where $b(3) > b(2)$. The $\beta$ in the time guarantees that if we work backwards from time $t/\ep_n + s_n$ to $t/\ep_n$, the probability a random walk particle will move by an amount much larger than $n^{b(2)/3}$, the size of the boxes at time $t/\ep_n$, is stretched exponentially small. See Lemma \ref{tailbound}. As in \cite{DurNeu} and \cite{CDP} this implies the conditional distribution of the position given that the lineage ends in a specific box is almost uniform, and hence the probability it lands on a 1 will be close to $u$.  A second consequence is that

\begin{lemma} \label{nuu}
With very high probability, the empirical finite dimension distributions at time $t/\ep_n + s_n$ will be close to $\nu_u(G_{x,y,v})$. 
\end{lemma}

\begin{proof}
To see this, note that we compute the probabilities of finite dimensional sets in the voter model equilibrium $\nu_u$ by starting the CRW with points at $y_0, \ldots y_m$, and running time to $s_n$. The particles that coalesce are a partition of the original set. We then flip a coin with a probability $u$ of heads (state 1)
to determine the states. Here we are only running time to $s_n$ so our partition is finer, but the final particles are roughly independent and uniform on the torus so whether they land on 1 or 0 are roughly independent coin flips. \end{proof}

The last paragraph shows that probabilities of the f.d.d.'s are close to the voter model equilibrium $\nu_u$. This enables us to conclude that the expected value of the drift of our process when the density is $x$ is close to $b(x)$. The next step is control the fluctuations about the mean. Using normal tail bounds on random walks in Lemma \ref{tailbound}, it follows that if $B_n$ is the event that some coalescing random walk at time $t/\ep_n+s_n$ moves by more than $n^{b(3)/3}$ in time $s_n$, then for any $\gamma>0$ we have for large $n$
\begin{align}
P(B_n) & \le n\exp( - (1-\gamma)  n^{2b(3)/3} / 2 n^{(2+\beta)b(2)/3} ) 
 \nonumber \\
& = n \exp\left( -\frac{1-\gamma}{2} n^{ [2b(3)-(2+\beta) b(2)]/3} \right)
\label{bdmove}
\end{align}

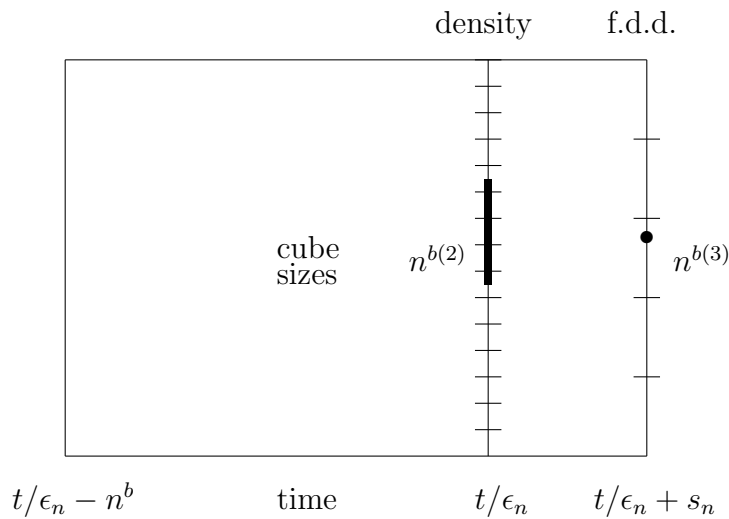
\begin{figure}[ht]
\begin{center}
\begin{picture}(300,230)
\put(40,40){\line(1,0){220}}
\put(40,190){\line(1,0){220}}
\put(40,40){\line(0,1){150}}
\put(200,40){\line(0,1){150}}
\put(260,40){\line(0,1){150}}
\put(120,115){cube}
\put(120,105){sizes}
\put(170,110){$n^{b(2)}$}
\put(270,110){$n^{b(3)}$}
\put(20,20){$t/\ep_n - n^b$}
\put(120,20){time}
\put(195,20){$t/\ep_n $}
\put(240,20){$t/\ep_n + s_n$}
\put(195,50){\line(1,0){10}}
\put(195,60){\line(1,0){10}}
\put(195,70){\line(1,0){10}}
\put(195,80){\line(1,0){10}}
\put(195,90){\line(1,0){10}}
\put(195,100){\line(1,0){10}}
\put(195,110){\line(1,0){10}}
\put(195,120){\line(1,0){10}}
\put(195,130){\line(1,0){10}}
\put(195,140){\line(1,0){10}}
\put(195,150){\line(1,0){10}}
\put(195,160){\line(1,0){10}}
\put(195,170){\line(1,0){10}}
\put(195,180){\line(1,0){10}}
\put(195,190){\line(1,0){10}}
\put(255,70){\line(1,0){10}}
\put(255,100){\line(1,0){10}}
\put(255,130){\line(1,0){10}}
\put(255,160){\line(1,0){10}}
\put(257,120){$\bullet$}
\put(180,200){density}
\put(245,200){f.d.d.}
\thicklines
\linethickness{1mm}
\put(200,105){\line(0,1){40}}
\end{picture}
\caption{Picture summarizing the proof. Here $s_n = n^{(2+\beta)b(2)/3}$.
The words at the top indicate the quantity that is ``good'' at each time, i.e., close to its average value on the cubes. The dark line at time $t/\ep_n$
shows the interval in which we will with high probability find the lineage of the black dot when it is worked backwards in time.}
\end{center}
\end{figure}

For the last inequality to be useful we need to choose $\beta$ so that $2b(3) - (2 + \beta) b(2) > 0$. The estimate in \eqref{bdmove} implies that the states of sites in cubes in the decomposition at time $t/\ep_n+ s_n$ that do not touch are independent on $B_n^c$.
We can divide our collection of cubes into 27 subcollections ${\cal C}_i$ of size $n^{1-b(3)}/27$ so that no two cubes in the subcollection touch. For $1\le i \le 27$, let $N_i$ be the number of times $G_{x,y,v}$ occurs in the union of the cubes in ${\cal C}_i$, let $N_{i,j}$ be the number of times $G_{x,y,v}$ occurs for $x$ in the $j$th cube in ${\cal C}_i$. If $x$ is close to the edge of the cube then some of the $x+y_i$ may be outside. However,	 the $y_i$ are fixed, so for large $n$ they will at worst be in an adjacent cube. 

For fixed $i$, the $N_{i,j}$ are independent on the event $B_n^c$, and $0 \le N_{i,j}/n^{b(3)} \le 1$.  Let $\rho_{i,j} =  EN_{i,j}/n^{b(3)}$. Let
$$
X_{i,j} = \frac{N_{i,j}}{n^{b(3)}} - \rho_{i,j} \in [-\rho_{i,j},1-\rho_{i,j}].
$$
Finally, let $\psi_{i,j}(\theta) = E\exp(\theta X_{i,j})$, let $Y_i = \sum_j X_{i,j}$, and let $M=n^{1-b(3)}/27$ be the number of cubes in each collection ${\cal C}_i$. If $\theta>0$, then, assuming $B_n^c$, we have
$$
e^{\theta M\eta} P( Y_i \ge M\eta ) \le \prod_{j} \psi_{i,j}(\theta),
$$
using the independence of the $N_{i, j}$ across $j$. So, we have
\begin{align}
P( Y_i \ge M\eta) &\le  e^{-\theta M\eta} \prod_{j} \psi_{i,j}(\theta)
\nonumber\\
&= \exp\left( M \left[ -\theta\eta + M^{-1}\sum_j \log \psi_{i,j}(\theta) \right] \right)
\label{ldfdd0}
\end{align}

Since we do not know much about $\psi_{i,j}(\theta)$, we will let $\eta_n = n^{-\alpha}$, and later choose $\theta_n$ so that $\lim_{n \to \infty} \theta_n = 0$. Expanding $\log\psi_{i,j}$ around 0:
\begin{align*}
\frac{d}{d\theta} \log \psi_{i,j}(\theta) &= \frac{\psi_{i,j}'(\theta)}{\psi_{i,j}(\theta)}, \\
\frac{d^2}{d\theta^2} \log \psi_{i,j}(\theta) &= \frac{\psi_{i,j}''(\theta)}{\psi_{i,j}(\theta)} - \frac{(\psi_{i,j}'(\theta))^2}{\psi_{i,j}^2(\theta)}. 
\end{align*}
When $\theta = 0$, we have $\psi_{i,j}(0)=1$ by definition, and also
\begin{align*}
\frac{d}{d\theta} \log \psi_{i,j}(0) &= EX_{i,j} = 0, \\
\frac{d^2}{d\theta^2} \log \psi_{i,j}(0) &= EX_{i,j}^2. 
\end{align*}
So, if $\theta_{i,n} \to 0$, then we have the approximation 
$$
\log \psi_{i,j}(\theta_{i,n}) \sim \frac{\theta_{i,n}^2}{2} EX^2_{i,j}.
$$
Since $X_{i,j} \in [-\rho_{i,j},1-\rho_{i,j}]$ and $EX_{i,j} = 0$,
$$
EX_{i,j}^2 \le  \rho_{i,j}(1-\rho_{i,j})
$$
To optimize the bound in \eqref{ldfdd0} we $d/d\theta$ the term in square brackets in \eqref{ldfdd0} to get
\beq
0=-\eta _n+ \theta_{i,n} M^{-1} \sum_j \rho_{i,j}(1-\rho_{i,j}), 
\label{chtheta}
\eeq
which says we want to take $\theta_n = \eta_n/\tau_i$,
where $\tau_{i} =  M^{-1} \sum_j \rho_{i,j}(1-\rho_{i,j})$. This gives the following large deviations bound 
\begin{align*}
P( Y_i \ge M\eta_n) &\le \exp\left( M \left[ -\frac{\eta_n^2}{\tau_i} + \frac{\eta_n^2}{2\tau_i^2}  \tau_i \right] \right) \\
& = \exp\left( - \frac{M \eta_n^2}{2\tau_i} \right) \le \exp( - M \eta_n^2 ),
\end{align*}
since $2\tau_i \le 1$.
The same reasoning can be used to get a bound on the other deviation. Since we have expanded the moment generating function around 
0 the bound is the same, giving the final result 
$$
P( |Y_i| \ge M\eta_n) \le \exp( - M \eta_n^2 )
$$
Define $Y = \sum_{i = 1}^{27} Y_i$, and then use the triangle inequality to get
$$
P( |Y| \ge 27 M\eta_n) \le 27\exp( - M \eta_n^2 )
$$
The last task is to relate this to the difference of the drifts. To do this, we note that
$$
Y = n^{-b(3)} \sum_{i,j} N_{i,j} - \sum_{i,j} \rho_{i,j}
$$
so we have 
$$
\frac{Y}{n^{1-b(3)}} = n^{-1} \sum_{x} 1(G_{x,y,v}) - \frac{1}{n^{1-b(3)}} \sum_{i,j} \rho_{i,j}
$$
Let $p^n_{x,y,v}$ be the probability of $G_{x,y,v}$ when we work backwards in the coalescing random walk starting from $x, x+y_1, \ldots x+y_k$ then we have
$$
\frac{1}{n^{1-b(3)}} \sum_{i,j} \rho_{i,j} = \frac{1}{n} \sum_x p^n_{x,y,v}
$$

In the three neighbor case we only have to consider: $y_1=e_1$, $y_2=e_2$, and $y_3=e_3$. When there are more neighbors, we have to consider a number of other possibilities, see the calculations in Section \ref{sec:cpert}. Let $r(v) = r(v_0,v_1,v_2,v_3)$ be the jump rate of vertex $x$ when the states are $v_i$. Multiplying by $r(v)$, summing over the relevant values of $y,v$ we have
$$
n^{-1} \sum_{x,y,v} 1(G_{x,y,v}) r(v) = \beta(\xi_{t/\ep_n + s_n})
$$
so we have
\beq
P\left( \left| \beta(\xi_{t/\ep_n + s_n}) - \frac{1}{n} \sum_{x,y,v} p^n_{x,y,v} r(v) \right| \ge 16 n^{-\alpha} \right) 
\le 27 \exp\left( - n^{1-b(3)-2\alpha}/27 \right)
\label{betatob}
\eeq
The choice of $s_n$ guarantees that as we work backwards in time the particles in the CRW move by an amount $\gg n^{b(3)}$. The bound in \eqref{mombd3} implies that each particle in the CRW lands on a 1 with probability close to $u$. It follows that  
$$
\left| \frac{1}{n} \sum_{y,v} p^n_{x,y,v} -b(u) \right| \le \eta/2
$$ 
with very high probability. The bounds derived above only works for fixed $t$.  However, it is easy to extend them so that they hold uniformly on $[0,t_0]$ 
and hence are valid for the integral. To do this, we subdivide the interval into subintervals of length $1/n^{1/2}\ep_n$. Within each interval the probability
there are more than $2n^{1/2}$ flips is $\le \exp(-c\sqrt{n})$. If we add this to previous error probability and multiply by the number of subinterval
we still have a result that holds with very high probability.

\subsection{Final details} \label{ss:fd}

To get long time survival, we will iterate. Let
$$
T_0 = \inf\{ t: |x_t - 1/2| < \eta \}
$$
and note that $x_t$ is the solution of the ODE so this is not random. 
Theorem \ref{DarNor} implies that  $|X(T_0)-1/2| \le 2\eta$ with very high probability. Let 
$$
T_1 = \inf\{ t > T_0 : |X_t - 1/2| \ge 4\eta \}
$$ 
and note that on $[T_0,T_1]$ we have $|X_t - 1/2| \le 4\eta$.  
There is a constant $t_\eta$ so that if $x(0)=1/2 + 4\eta$ or $x(0)=1/2 - 4\eta$ then
$|x(t_\eta)-1/2| \le \eta$. Let $S_1=T_1+t_\eta$. Since $T_1$ is random, $S_1$ is a random time. However, due to the Markov process, we can translate time to apply Theorem \ref{DarNor} again. That is, consider $\tilde{X}_t := X_{t + T_1}$. Then since $| \tilde{X}_0 - 1/2 | = 4 \eta$, Theorem \ref{DarNor} implies that with high probability $|\tilde{X}_{t_\eta} - 1/2| = |X(S_1)-1/2| \le 2\eta$ and
$|X_t -1/2| \le 5\eta$ on $[T_1,S_1]$. For $m \ge 2$, let 
$$
T_m = \inf\{ t > S_{m-1} : |X_t - 1/2| \ge 4\eta \}\quad\hbox{and}\quad S_m=T_m+t_\eta. 
$$
We can with high probability iterate the construction $n^{k}$ times before it fails. Since each cycle takes at least $t_0$ units of time, taking $\eta=\gamma/5$ the proof of Theorem \ref{persist} is complete.

\clearp

\section{Rapid extinction for $q>1$} \label{sec:rapext}

In this section we will prove Theorem \ref{dieout}. There are two steps to the proof. First, we use the results in Section 4 to show that the fraction of 1'sin the random process is close to solution of the ODE until time
\beq
\tau = \min\{ t: x_t < n^{-(1-b(0)} \},
\label{taudef}
\eeq
where $b(0)$ will be defined in the proof of Lemma \ref{lem:dyingFirst}. The second step is to prove that when we start with $\le n^{b(0)}$ ones, then fluctuations in the voter model will cause it to hit $0$ in time $\le C n^{b(0)}$. This time is $<n^b$ for large $n$, so by results in Section \ref{ss:igbr}, it is legitimate to assume that the process acts like the voter model. The proof for the second step is based on a Green's function calculation and estimates for the rate of change of the number of ones in the voter model.

\subsection{First step}

\begin{lemma} \label{lem:dyingFirst}
Suppose $X_0 < 1/2$ and let $\tau$ be defined in \eqref{taudef}. Then, for any $\eta > 0$, as $n \to \infty$,
\[
\P \left( |X_\tau - n^{-(1 - b(0))} | < \eta n^{-(1 - b(0))} \right) \to 1. 
\]
\end{lemma}

\begin{proof}
We use \eqref{mombd3} from Section \ref{ss:TvsT}. If $X_0=u$ and we divide the torus at time $t/\ep_n$ into boxes of side $r=n^{b(2)/3}$, then taking $m$ large in\eqref{mombd3} gives
\beq 
P\left( \hbox{ for some $i$ }\left| \frac{N_{r,i}}{n^{b(2)}} - u \right| \ge [u(1-u)]^{1/2} n^{-b(2)/6+\beta} \right)  \le n^{-k},
\eeq 
for any $\beta > 0$ and $k < \infty$. Since $u^{1/2} > (u(1 - u))^{1/2}$, we can change this to 
\beq 
P\left( \hbox{ for some $i$ }\left| \frac{N_{r,i}}{n^{b(2)}} - u \right| \ge u^{1/2} n^{-b(2)/6+\beta} \right)  \le n^{-k}.
\label{mombdv}
\eeq 
For this estimate to be useful, we need $u \gg u^{1/2} n^{-b(2)/6+\beta}$ which is equivalent to $u \gg n^{-b(2)/3 + 2\beta}$. If $b(2)$ is close to 1 and $\beta$ is small, we can define $b(0)$ by
$$
1 - b(0) = b(2)/3 - 2\beta,
$$
so that $b(0) < \min\{b,1-\alpha\}$ where $\alpha>1/3$ is the quantity from Theorem \ref{dieout}.  Combining these estimates and using results from the previous section we have that if $x_0 < 1/2$ and $\eta>0$ then as $n\to\infty$
$$
P( |X_t - x_t| \le \eta x_t \hbox{ for all $t\le \tau$}) \to 1 .
$$
Lemma \ref{lem:dyingFirst} follows. \end{proof} 

This result shows that the number of 1's gets driven to $\le (1+\ep) n^{b(0)}$ at the deterministic time $\tau$. To complete the process of extinction we will rely on fluctuations
in the voter model.

\subsection{Green's function calculation} \label{ss:greenf}

To motivate the calculation in the next lemma we note that the voter model is a time change of simple random walk.

\begin{lemma}
Let $S_t$ be continuous-time simple random walk on $\{0, \dots, n \}$ with jump-rate $r(j)$ at position $j$. Let $0 < x < z \le n$ be integers, and $T_{0, z}$ the first time that $S_t$ hits $0$ or $z$. Then,
\begin{equation} \label{gf}
E_x T_{0, z} =  \sum_{y = 1}^{x} \frac{2y}{r(y)} + \sum_{y = x+1}^{z} \frac{2x}{r(y)} 
- \sum_{y = 1}^{z} \frac{2xy}{z r(y)}.
\end{equation}

\end{lemma}

\noindent 
Since $P_x(T_z< T_0) = x/z$, this is enough to bound the extinction time if $x/z\to 0$.

\begin{proof}
First consider the embedded discrete-time chain of $S_t$. For $0 \le y \le z$, let $N_x(y)$ be the number of times the random walk visits $y$ before hitting $0$ or $z$, starting from position $x$. Consider the Green's function
\[
G_0(x, y) = \E[N_x(y)].
\]
Fix $y$ and write $g(x) = G_0(x, y)$. Then we have that $g$ satisfies
\[
\begin{cases}
g(0) = 0 \\
g(x) = \frac{1}{2} \left( g(x + 1) + g(x-1) \right), & x \neq 0, y, z \\
g(y) = 1 + \frac{1}{2} \left( g(y + 1) + g(y-1) \right) \\
g(z) = 0
\end{cases}.
\]
From this it is clear that $g$ should be linear and increasing on $[0, y]$ and linear and decreasing on $[y, z]$. That is,
\[
\begin{cases}
g(x) = c_1 x  & 0 \le x \le y \\
g(x) = c_2(z - x) & y \le x \le z.
\end{cases}.
\]
To satisfy the conditions for $g(x)$ and $g(y)$, the constants must be 
\begin{align*}
c_1 = \frac{2(z - y)}{z}, \,\,\,\,\,
c_2 = \frac{2y}{z}.
\end{align*}

The walk will spend an average of $1/r(y)$ units of time at position $y$ before jumping. Thus, if $G(x, y)$ is defined to be the expected amount of time the continuous time walk spends at $y$, started from $x$, before hitting $0$ or $z$, we have:
\[
G(x, y) = \frac{1}{r(y)} \cdot G_0(x, y)  = \frac{1}{r(y) } \cdot \begin{cases}
2x(z - y)/z  & x \le y \\
2(z-x)y/z  & x \ge y
\end{cases}
\]
Thus, the expected total time before being absorbed, started from $x$, is
\begin{align*}
E_x[T_{0, z}] &= \sum_{y = 1}^z G(x, y) = \sum_{y = 1}^{x} \frac{2y}{z} \cdot (z - x) \cdot \frac{1}{r(y)} + \sum_{y = x+1}^z \frac{2(z - y)}{z} \cdot x \cdot \frac{1}{r(y)} \\
&= \sum_{y = 1}^{x} \frac{2y}{r(y)} + \sum_{y = x+1}^z \frac{2x}{r(y)} - \sum_{y = 0}^z \frac{2xy}{z r(y)},
 \end{align*}
which establishes \eqref{gf}
\end{proof}

\subsection{Boundary size calculations} \label{ss:boundary}

To use \eqref{gf} to bound the extinction time, we need to understand the size of the boundary of the voter model: 
$\partial\xi = \{ \{x,y\} : x \sim y, \, \xi(x) \neq \xi(y) \}$.
Here $x\sim y$ means that $x$ and $y$ are neighbors and $\{x,y\}$ is the un-oriented edge that connects them. For a voter model configuration $\xi$, let $|\xi| = \sum_x \xi(x)$ be the number of 1s. The next result gives trivial upper and lower bounds on $|\partial\xi|$ when  $|\xi| = k$:
\begin{equation} 
\label{easybds}
C_d k^{1/d} \le | \partial \xi| \le 2d k. 
\end{equation}

Using \eqref{gf}, we see that if $x=n^p$ and $z=n^q$ for some $0 < p < q < 1$, then for $r(y) = y$,
\beq
 E_x T_{0, z} \le \sum_{y = 1}^{x} \frac{2y}{y} + \sum_{y = x+1}^{z} \frac{2x}{y}
 \le C x + 2x [\log(z)-\log(x)] \le C'x \log(z)
\label{ub0}
\eeq
If $p=b(0)$ and $q> p$, this gives us what we want, an extinction time $\ll n^b$.

On the other hand, if we use the lower bound and plug in $r(y) = y^{1/3}$, then
\beq
E_x T_{0, z} \le \sum_{y = 1}^{x} \frac{2y}{y^{1/3}} + \sum_{y = x+1}^{z} \frac{2x}{y^{1/3}}
\le C ( x^{5/3} + x^{2/3}z) \le C' x^{2/3} z
\label{ub1}
\eeq
If we take $x=n^{b(0)}$ and $z=n^c$ then this is $\le C n^{5b(0)/3}$, which is much longer than the interval of length $n^{b}$ over which the process behaves like the voter model. Combining \eqref{easybds} and \eqref{ub1} gives

\begin{lemma} \label{smallk}
If $x=n^p$ with $p<3b/5$ and $z=n^q$ with $q>p$ and $2p/3+ q < b$ then
$$
P_x( T_{0,z} \le n^b ) \to 1 \quad\hbox{as $n\to\infty$}.
$$ 
\end{lemma}  

\noindent
This will let us show that the time spent at small values of $|\partial \xi_t|$ can be ignored. For larger values, we 
need a more precise statement about the size of the boundary. This has been done by Cox, Durrett, and Perkins \cite{CDP2},
in order to show that in $d\ge 2$ the rescaled voter model converged in distribution to super-Brownian motion. 
This was later used by  Bramson, Cox, and LeGall \cite{BCL} to prove a result for the voter model in $d\ge 3$ started at 0. 
See Theorem 4 on page 1012 in \cite{BCL}. 

To prepare for stating our lemma we describe the result from \cite{CDP2}. 
They use a general probability kernel $p(z)$. In our case $p(z)=1/6$ for the nearest neighbors of 0. 
If $\xi_t(x)=1$ we let
$$
V_t(x) = \sum_y p(y-x) 1_{(\xi_t(y) = 1)}
$$
If $\xi_t(x)=0$ we set $V_t(x)=0$. This part of the definition is not really needed in the statement since $X^N_s$ is supported by points
on the rescale lattice in state 1. On page 202 of their result you find the following result.

\mn
(I1) There is a finite $\gamma>0$ so that for all $\phi \in C_0^\infty(\RR^d)$ and $T>0$
$$
E\left[ \left( \int_0^T X^N_s( [V_{N,s}-\gamma] \phi^2 ) \, ds \right)^{\sqz 2} \right] \to 0
$$
Here $X^N_t$ is the voter model with space scaled by $\sqrt{N}$ and time scaled by $N$ and turned into a measure by assigning mass $1/N$ to states
in state 1, see (1.4), and $V_{N,s}(x)$ is a suitably rescaled version of $V_t(x)$. The formula on page 202 has $V'$ because they want to write
the formula so that it is valid for $d=2$ and $d \ge 3$. 

In our situation $\gamma = 2d\beta_d$. However, in this proof we need control on the size of the error. The reader should think of $s$ as a point in the time
interval $[t/\ep_n - n^b/2, t/\ep_n]$ over which our process behaves like the voter model. 

\begin{lemma}\label{bdyerr}
If $k$ is large and the density of 1's is small then
$$
P\left( \left. \frac{|\partial \xi_s|}{ | \xi_s|} \not\in [(1-\ep) 2d\beta_d k , (1+\ep) 2d\beta_d k]  
\, \right| \, |\xi_s| = k \right) \le k^{-2/3}
$$
\end{lemma}

\begin{proof} Pick a site $x$ at time $s$ with $\xi_s(x)=1$. When this holds the coalescing random walk starting at $x$ at time $s$ lands on a site in state $1$ at time $t/\ep_n - n^b$. Let $r=k^\alpha$ where $\alpha$ is small and follow the CRW path backwards in time for $r$ units of time. If we let $h(s,x)$ be the probability the CRW starting at time $s$ lands on 1 at time $t/\ep_n - n^b$, then an elementary conditional probability shows that the probability our conditioned CRW particle at $x$ at time $s$ is at $y$ at time $s-r$ is
$$
\bar p_{s,r}(x,y) = p_r(x,y) \frac{h(r,y)}{h(s,x)}
$$
This result is often known as Doob's $h$-transform. Since the lineage will wrap around the torus in the remaining $\ge n^b/2$ units of time, the ratio is close to 1 and can be ignored. 

For each neighbor $y$ of an $x$ with $\xi_t(x)=1$, let $V_{x,y}=1$ if it does not coalesce with $x$ by time $r$ and 0 otherwise. For any $\alpha>0$, if $k$ is large and the density of 1's is $u$ which is small then  
$$
\left| \frac{P(V_{x,y}=1)}{ 2d\beta_d k} - 1\right| < \eta/2.
$$  
Here we are using the hydrodynamic limit Lemma \ref{nuu} to conclude that the distribution of the process is close to $\nu_u$ at time $r$. 

Let $W_x = \sum_{y\sim x} V_{x,y}$,  $\mu(x) =  \sum_{y\sim x} EV_{x,y}$, and 
$$
S_k = \sum_x^\star \bar W_x \quad\hbox{where} \quad W_x - \mu(x)
$$
where $\Sigma^\star_x$ is short for $\sum_{x: \xi^0_t(x)=1}$. Arguments in Section \ref{ss:diffm} imply that if $|x-x'| > s$ then the correlation
between $W_x$ and $W_{x'}$ is small enough to be ignored so
$$
E(S_k^2)  = \sum^\star_x \sum^\star_y E[\bar W_x \bar W_y]  \le 36 k \cdot Cr^3
$$
since $|\bar W_x | \le 6$ and for a given $x$ there are at most $Cr^3$ values of $y$ with $|x-y| \le r$. If we use Chebyshev's inequality
$$
P( |S_k| \ge \ep k ) \le \frac{36 k \cdot Cs^3}{ k^2}   \le C k^{-1+3\delta}
$$
If $\alpha < 1/10$ this gives the desired result.
\end{proof}

\subsection{Extinction time}
\label{ss:exttime}

The results about the boundary of the voter model can now be applied to the Green's function calculation to get the result

\begin{lemma}  \label{lem:extinctionTime}
Consider the voter model started with configuration $|\xi_0| = x$ and let $T_{0, z}$ be the first time the configuration hits $0$ or $z$. If $x=n^{b(0)}$ and $z=n^c$ with $c>b(0)$  then 
$$
E_x[T_{0,z}] \le C n^{b(0)}
$$
\end{lemma}

\begin{proof} 

We can divide the sum in \eqref{gf} into the pieces where Lemma \ref{smallk} can be applied. That is, define $x' = n^{p} < x$ so that $p < 3 b(0)/5$ and $2p/3 + c < b(0)$. Then,
\[
\E_x[T_{0, z}] \le \E_{x'}[T_{0, z}] + \E_x[T_{x', z}].
\]
The first term is less than a constant times $n^{b(0)}$ by Lemma \ref{smallk}. To bound the second hitting time, we use \eqref{ub0} and Lemma \ref{bdyerr} to conclude that the expected amount of time
when $|\partial\xi_s|/|\xi_s|$ is not within $\ep$ of $2d\beta_d$ is
$$
\le  \sum_{y = 1}^{x} \frac{2y}{y^{1/3}} y^{-2/3} + \sum_{y = x+1}^{z} \frac{2x}{y^{1/3}} y^{-2/3} \le \sum_{y = 1}^{x} \frac{2y}{y^{1/3}}  + \sum_{y = x+1}^{z} \frac{2x}{y^{1/3}} \le  Cn^b
$$
which finally completes the proof.
\end{proof}

Theorem \ref{dieout} now immediately follows: apply Lemma \ref{lem:dyingFirst} to get that $U_n(\alpha \log n) < n^{-(1-b(0))}$ with high probability. Next, use Section \ref{ss:igbr} so that with high probability we can assume the $q$-voter model only experiences voter branching events for the remainder of the time. Lemma \ref{lem:extinctionTime} then proves that with high probability the unscaled voter model started with $n^{b(0)}$ occupied sites will hit $0$ or $n^c$ in an additional time of $C n^{b(0)}$. The probability that the process hits $0$ first is simply $(n^c - n^{b(0)})/n^c \to 1$. Since $b(0) > 2/3$, this additional time is $o(1)$ for the time-scaled process $U_n(t)$. Thus,
$$
P \left( U_n(\alpha \log n)=0 \right) \to 1 \qquad\hbox{as $n\to\infty$.}
$$

\clearp

\section{Computing the perturbation} \label{sec:cpert}

In this section, Theorem \ref{limitODE} is proved. Recall Theorem \ref{detlim} state that the limiting ODE for the model with a $k$-sized neighborhood is 
\[
\frac{du}{dt} = \sum_{m = 1}^{k-1} r_i^k (\rho_m^0(u) - \rho_m^1(u)),
\]
where $\rho_m^i(u)$ is the probability under the voter model equilibrium $\nu_u$ that the origin is in state  $i$ and a exactly $m$ of the neighbors are in state $1-i$. In this section, we analyze these quantities. Before stating the proof for a general $k$, we first show an explicit proof for a neighborhood of size $3$ to give a flavor of how the individual terms are computed, while introducing some necessary notations in an organic manner.

\subsection{\bf k=3}

To compute $\rho^0_i$ we have to compute the coalescence fate of 0, $e_1$, $e_2$, $e_3$. There are 7 possibilities
\begin{center}
\begin{tabular}{lcccc}
one & 0 ; 3 & 1: 2 & 2 ; 1  & 3; 0\\
two & 0; 2, 1 & 1: 1, 1  \\
three & 0; 1,1,1 
\end{tabular}
\end{center}

The first number in each string gives the number of neighbors that coalesce with 0. The others give the size of the limiting coalescing clusters formed by the
remaining neighbors. The word at the beginning of the row is the number of numbers after the semi-colon. We can ignore $3;0$ because in that case all the neighbors have the same state as 0. 

Let $\rho^0_i$ be the probability that in the voter equilibrium $\nu_u$ the origin is 0 while exactly $i$ of the neighbors are 1. Factoring out the probability the origin is we have $\rho^0_i = (1-u)q_i(u)$.To compute the $q_i(u)$ we use the following table. 

\begin{itemize}
  \item The coefficients of $u$ come from the ``one'' terms.
  \item The coefficients of $u^2$ and $u(1-u)$ come from the ``two'' terms. There is no $(1-u)^k$ since all the neighbors would be 0. $p(1;1,1)$ appears three times since only 0,0 is impossible. $p(0;2,1)$ only appears twice since 0,0 and 1,1 are impossible.
  \item The coefficients of $u^2(1-u)$ and $u(1-u)^2$ come from the ``three'' terms. There is no $u^3$ or $(1-u)^3$ since all neighbors would be 0 or 1. For this reason $p_{0;1,1,1}$ appears $2^3 - 2 = 6$ times
  \end{itemize} 

The meaning of the first column will become clear when the reader reaches \eqref{diffq}
\beq
\begin{matrix}
\Delta_i(u) & term  & q_1(u) & q_2(u) \\
0 & u & p_{2;1} & p_{1;2} \\
-1 & u^2 &  & p_{1;1,1}  \\
1 & u(1-u) & p_{0;2,1} + 2p_{1;1,1} & p_{0;2,1}\\
1 & u(1-u)^2 & 3p_{0;1,1,1} \\
0 & u^2(1-u) & & 3p_{0;1,1,1} 
\end{matrix}
\label{table3}
\eeq
so reading down the columns we have
\begin{align*}
q_1(u) & = p_{2,1} u + [2p_{1,1,1}+ p_{0,2,1}] u(1-u) + 3p_{0,1,1,1}u(1-u)^2 \\
q_2(u) & = p_{1,2} u +p_{1,1,1}u^2 + p_{0,2,1} u(1-u) + 3p_{0,1,1,1} u^2(1-u) 
\end{align*}

Let $\rho^1_i$ be the probability that in the voter equilibrium $\nu_u$ the origin is 1 while exactly $i$ of the neighbors are 0. From the previous calculation we see that $\rho^1_i = u q_i(1-u)$ so we have
$$
\langle h_{0,1} - h_{1,0} \rangle_u = \sum_{i=1}^2 r_i (\rho^0_i - \rho^1_i)
$$
The quantity in parentheses is $\Delta_i(u) \equiv (1-u)q_i(u) - u q_i(1-u)$. 
Taking difference we have (the first column indicates the term in $q_i(u)$)
\begin{align}
u \quad & \quad u(1-u) - (1-u) u  = 0 
\nonumber\\
u^2 \quad & \quad u^2(1-u) - (1-u)u^2 = u(1-u)(2u-1)
\nonumber\\
u(1-u) & \quad u(1-u)^2 - u^2(1-u)  = u(1-u)(1-2u) 
\label{diffq}\\
u(1-u)^2 & \quad u(1-u)^3 - u^3(1-u)  = u(1-u)[ (1-u)^2 - u^2 ] = u(1-u)(1-2u)
\nonumber\\
u^2 (1-u) & \quad u^2(1-u)^2 - (1-u)^2u^2  = 0 \nonumber
\end{align}
so consulting \eqref{table3} we have
\begin{align*}
\Delta_1(u) & = [2p_{1,1,1}+ p_{0,2,1} + 3 p_{0,1,1,1}] u(1-u)(1-2u) \\
\Delta_2(u) & = [- p_{1,1,1} + p_{0,2,1}] u(1-u)(1-2u)
\end{align*}
and the reaction term is 
\begin{align*} 
\frac{\phi(u)}{u(1-u)(1-2u)}  & = r_1[2p_{1,1,1} + p_{0,2,1} + 3 p_{0,1,1,1} ] \\
& + r_2 [ - p_{1,1,1} + p_{0,2,1} ] 
\end{align*}
If $2r^3_1>r^3_2$ so the right-hand side is positive. Using \eqref{rforqv} we see that in the q-voter model with $q<1$
$$
2r^3_1 = 2/3 \log (3) > 2/3 \log(3/2) = r^3_2
$$
so the reaction term is $c_3 u(1-u)(1-2u)$ with $c_3>0$. When $q>1$ the reaction term is $-c_3 u(1-u)(1-2u)$.

\subsection{\bf General k}

In this case we have to compute the coalescence fate of $0$ with $k$ neighbors. Again $\rho^0_i=(1-u)q_i(u)$, where the functions  $q_i(u)$, $i\leq k-1$  defined as before are   polynomials with terms of the type $u^a(1-u)^b$. First let us look at the difference $\Delta_{a,b} (u)$ of these terms, where $\Delta_{a,b}(u)=\rho_i^0-\rho_i^1=u^a(1-u)^{b+1}-u^{b+1}(1-u)^a$. Note that $\Delta_{a,b} (u)=0$ if $a=b+1$. 

In the case $a\leq b$ we have
\begin{align*}
\Delta_{a,b} (u)&= u^a(1-u)^{b+1}-u^{b+1}(1-u)^a\\
&=u^a(1-u)^a[(1-u)^{b-a+1}-u^{b-a+1}]\\
&=u^a(1-u)^a(1-2u)\left[\sum_{j=0}^{b-a}u^j(1-u)^{b-a-j}\right].
\end{align*}
To see the last step write $1-2u = (1- u)-u$ and the telescope the sum. In the case $a>b+1$
\begin{align*}
\Delta_{a,b} (u)&= u^a(1-u)^{b+1}-u^{b+1}(1-u)^a\\
&=u^{b+1}(1-u)^{b+1}[u^{a-b-1}-(1-u)^{a-b-1}]\\
&=-u^{b+1}(1-u)^{b+1}(1-2u)\left[\sum_{j=0}^{a-b-2}u^j(1-u)^{a-b-2-j}\right]
\end{align*}
Since $\sum_{j=0}^n u^j (1-u)^{n-j}>0$ on $[0,1]$ we have that $0,1$ and $1/2$ are the only roots of $\Delta_{a,b} (u)$. Also note that $\Delta_{a,b} (u)=-\Delta_{b+1,a-1} (u)$. We claim  
$$
\frac{\phi(u)}{u(1-u)(1-2u)}=f(u), 
$$
where $f(u)$ is a positive polynomial  in $u$ with no real roots. To prove this, given a coalescence fate $s_0; s_1,s_2,s_3,\cdots, s_j$ where $\sum_j s_j=k$ we look at number of ways to obtain $a$ clusters with opinion $1$ (which gives the coefficients of the terms  $u^a(1-u)^b$, $a>b+1$) and compare it  with the  number of ways to obtain $b+1$  clusters with opinion $1$ (which gives the coefficients of the terms $u^{b+1}(1-u)^{a-1}$). 

First, suppose $b=0$ and $a\geq 2$.  Let $s_0$ be the number of neighbors that have coalesced with $0$, and $s_1,s_2,\cdots, s_a$ be the sizes of the limiting coalescing clusters formed by the rest of the neighbors, where we assume that the sizes are arranged in an increasing order, i.e., $s_1\leq s_2\leq\cdots\leq s_a$.  The coefficient of  $\Delta_{a,0} (u)$ in $\phi(u)$ is given by $r_{s_1+\cdots+s_a}p_{s_0;s_1,\cdots,s_a}$(Since all the clusters have opinion 1,  there is only one way to choose). Similarly the coefficient of  $\Delta_{1,a-1} (u)$ in $\phi(u)$ is given by $(r_{s_1}+\cdots+r_{s_a})p_{s_0;s_1,\cdots,s_a}$ (Since exactly one of the clusters has opinion 1, there are $a$ different choices, the coefficient of each of the clusters needs to be added individually). 
 
 Since $s_i$'s are increasing in $i$, so
 \begin{align*}
 \log(k/s_a)\leq\log (k/s_j)\qquad \forall j\in\{1,2,\cdots, a-1\}.
 \end{align*}
 So by the definition $r_i^k=\frac{i}{k}\log(k/i)$, and using the inequality above we have
 \begin{align*}
 r_{s_1+\cdots+s_a}&=\frac{s_1+\cdots+s_a}{k}\log(k/{(s_1+\cdots+s_a)})\\
 &\leq \frac{s_1+\cdots+s_a}{k}\log(k/{s_a})\\
 &=\frac{s_1}{k}\log(k/{s_a})+\frac{s_2}{k}\log(k/{s_a})+\cdots+\frac{s_a}{k}\log(k/{s_a})\\
 &\leq \frac{s_1}{k}\log(k/{s_1})+\frac{s_2}{k}\log(k/{s_2})+\cdots+\frac{s_a}{k}\log(k/{s_a})\\
 &=r_{s_1}+r_{s_2}+\cdots+r_{s_a}.
 \end{align*}
 Since $\Delta_{a,0} (u)=-\Delta_{1,a-1} (u)$, if we only look at terms of the type $\Delta_{1,a-1} (u)p_{s_0;s_1,\cdots,s_a}$ (which is non-negative) in $\phi(u)$, we get a non-negative polynomial in $u$ with no roots other than $0,1$ and $1/2$.
 
 Now suppose $b\neq 0$ and $a\geq b+2$. As explained in the previous case, let $s_0$ be the number of neighbors that coalesce with $0$, and  $s_1,s_2,\cdots, s_{a+b}$ be the sizes of the limiting coalescing clusters formed by the rest of the neighbors, where we assume that the sizes are arranged in an increasing order, i.e., $s_1\leq s_2\leq\cdots\leq s_{a+b}$. There are $\binom{a+b}{a}$ ways of choosing $a$ clusters out of the $a+b$ clusters. Denote the total size of each of these clusters by $x_i$, where $1\leq i\leq \binom{a+b}{a}$, where wlog we assume that the sizes are arranged in an ascending order. The coefficient of  $\Delta_{a,b} (u)$ in $\phi(u)$ is given by $p_{s_0;s_1,s_2,\cdots,s_{a+b}}\sum_{i=1}^{\binom{a+b}{a}} r_{x_i}$. Given $1\leq i\leq a+b$, the number of clusters in which cluster $s_i$ has opinion $1$ is given by $\binom{a+b-1}{a-1}$. Hence the total size of all the clusters, where $a$ of them have opinion $1$, is given by 
 $$\sum_{i=1}^{\binom{a+b}{a}}x_i=\binom{a+b-1}{a-1}\left(s_1+s_2+\cdots + s_{a+b}\right).$$
 
 Using a similar argument there are $\binom{a+b}{b+1}$ ways of choosing $b+1$ clusters out of the $a+b$ clusters. Denote the total size of each of these clusters by $y_i$, where $1\leq i\leq \binom{a+b}{b+1}$, where wlog we assume that the sizes are arranged in an ascending order. The coefficient of  $\Delta_{b+1,a-1} (u)$ in $\phi(u)$ is given by $p_{s_0;s_1,s_2,\cdots,s_{a+b}}\sum_{i=1}^{\binom{a+b}{b+1}} r_{y_i}$. Given $1\leq i\leq a+b$, the number of clusters in which cluster $s_i$ has opinion $1$ is given by $\binom{a+b-1}{b}=\binom{a+b-1}{a-1}$. Hence the total size of all the clusters, where $b+1$ of them have opinion $1$, is given by 
 $$\sum_{i=1}^{\binom{a+b}{b+1}}y_i=\binom{a+b-1}{a-1}\left(s_1+s_2+\cdots + s_{a+b}\right).$$
For ease of notation, let us denote $\binom{a+b}{a}$ by $n$ and $\binom{a+b}{b+1}$ by $m$.
Then $m>n$ since
\begin{align*}
\binom{a+b}{b+1}-\binom{a+b}{a}=&\binom{a+b}{a}\left(\frac{a}{b+1}-1\right)\\
=&\binom{a+b}{a}\left(\frac{a-b-1}{b+1}\right)>0.
\end{align*}

 Since $\sum_{i=1}^{n}x_i=\sum_{i=1}^{m}y_i$, and the $x_i$s as well as the $y_i$ s are arranged in ascending order, we have $x_i>y_i+m-n$, for $1\leq i\leq n$.
 \begin{align*}
 &x_1\log x_1+ x_2\log x_2+\cdots +x_n\log x_n -y_1\log y_1-y_2\log y_2-\cdots-y_m\log y_m\\
 >&x_1\log x_1+ x_2\log x_2+\cdots +x_n\log x_n -y_1\log y_1-y_2\log y_2-\cdots-y_m\log x_n\\
 =&x_1\log x_1+ x_2\log x_2+\cdots +x_n\log x_n -y_1\log y_1-y_2\log y_2-\cdots-y_{j-1}\log y_{j-1}-c\log y_j,
 \end{align*}
 where  $y_n+y_{n+1}+\cdots +y_j -c=x_n$. Now we have $\sum_{i=1}^{n-1}x_i=c+\sum_{i=1}^{j-1}y_i$. Repeating the same process as explained above $n-1$ times, we have 
 $$\sum_{i=1}^{n}x_i\log x_i>\sum_{i=1}^{m}y_i\log y_i.$$
 Now using the definition of $r_i^k$
 \begin{align*}
 \sum_{i=1}^{n}r^k_{x_i}=&\sum_{i=1}^{n}\frac{x_i}{k}\log(k/x_i)
 =\sum_{i=1}^{n}\frac{x_i}{k}\left[\log(k) -\log (x_i)\right]\\
 =&\sum_{i=1}^{m}\frac{y_i}{k}\log(k)-\sum_{i=1}^{n}\frac{x_i}{k}\log(x_i)
 <\sum_{i=1}^{m}\frac{y_i}{k}\log(k)-\sum_{i=1}^{m}\frac{y_i}{k}\log(y_i)\\
 =&\sum_{i=1}^{m}\frac{y_i}{k}\log(k/y_i)
 = \sum_{i=1}^{m}r^k_{y_i}.
 \end{align*}
 Now using the above inequality along with the fact that $\Delta_{a,b}=-\Delta_{b+1,a-1}$ , if we only look at terms of the type $\Delta_{b+1,a-1} (u)p_{s_0;s_1,\cdots,s_{a+b}}$ (which is non-negative) in $\phi(u)$, we get a non-negative polynomial in $u$ with no roots other than $0,1$ and $1/2$. This proves Theorem \ref{limitODE} for $q < 1$. 
\begin{corollary}
Fix $q>1$. For a $q$-voter model with $k$-neighbors, the reaction function defined in \eqref{ODErhs} simplifies to 
\beq
\phi(u)=-u(1-u)(1-2u)f^k(u),
\eeq
where $f^k(u)$ is a strictly positive polynomial in $u$.
\end{corollary}

\begin{proof}
Recalling the perturbation from \eqref{rforqv} and \eqref{rforq>1}, note that the perturbation when $q>1$ has the same value as the perturbation when $q< 1$ but with the opposite sign. This along with the above work proves the corollary.   
\end{proof}

\section*{Acknowledgments}
This work was begun during the 2019 AMS Math Research Communities meeting on Stochastic Spatial Models, June 9-15, 2019. We would like to thank
Hwai-Ray Tung, a graduate student at Duke for producing the figures. RD was partially supported by NSF grant DMS 1809967 from the probability program. MS was supported by a National Defense Science \& Engineering Graduate Fellowship. PA was partially supported by the NSF Grant DMS 1407504.

\clearp

\end{document}